\newif\ifmarek
\marektrue
\documentclass[reqno,twoside,11pt]{amsart}

\usepackage{srcltx}
\usepackage{amsmath}
\usepackage{amsfonts}
\usepackage{amssymb}
\usepackage{verbatim}
\usepackage{epsfig}

\IfFileExists{epsf.def}{\input epsf.def}{\usepackage{epsf}}

\DeclareFontFamily{OT1}{eusb}{} \DeclareFontShape{OT1}{eusb}{m}{n} {<5> <6> <7> <8> <9> <10> <11> <12> <14.4> eusb10}{}
\DeclareMathAlphabet{\eusb}{OT1}{eusb}{m}{n}

\DeclareFontFamily{OT1}{eusm}{} \DeclareFontShape{OT1}{eusm}{m}{n} {<5> <6> <7> <8> <9> <10> <11> <12> <14.4> eusm10}{}
\DeclareMathAlphabet{\eusm}{OT1}{eusm}{m}{n}

\DeclareFontFamily{OT1}{eufm}{} \DeclareFontShape{OT1}{eufm}{m}{n} {<5> <6> <7> <8> <9> <10> <11> <12> <14.4> eufm10}{}
\DeclareMathAlphabet{\mathfrak}{OT1}{eufm}{m}{n}

\DeclareFontFamily{OT1}{fraktura}{}
\DeclareFontShape{OT1}{fraktura}{m}{n} {<5> <6> <7> <8> <9> <10> <11> <12> <13> <14.4> [1.1] eufm10}{}
\DeclareMathAlphabet{\fraktura}{OT1}{fraktura}{m}{n}

\DeclareFontFamily{OT1}{cmfi}{} \DeclareFontShape{OT1}{cmfi}{m}{n} {<5> <6> <7> <8> <9> <10> <11> <12> <13> <14.4> [0.9] cmfi10}{}
\DeclareMathAlphabet{\cmfi}{OT1}{cmfi}{b}{n}

\DeclareFontFamily{OT1}{cmss}{} \DeclareFontShape{OT1}{cmss}{m}{n} {<5> <6> <7> <8> <9> <10> <11> <12> <13> <14.4> cmss10}{}
\DeclareMathAlphabet{\cmss}{OT1}{cmss}{m}{n}

\ifmarek
\IfFileExists{mathptmx.sty}{\usepackage{mathptmx}\usepackage{mathrsfs}}
{\usepackage{times}\usepackage{mathrsfs}}
\renewcommand{\mathcal}{\eusm}
\fi

\newtheoremstyle{thm}{1.5ex}{1.5ex}{\itshape\rmfamily}{} {\bfseries\rmfamily}{}{2ex}{}

\newtheoremstyle{def}{1.5ex}{1.5ex}{\slshape\rmfamily}{} {\bfseries\rmfamily}{}{2ex}{}

\newtheoremstyle{rem}{1.5ex}{1.5ex}{\rmfamily}{} {\bfseries\rmfamily}{} {1.5ex}{}

\newenvironment{proofsect}[1] {\vskip0.1cm\noindent{\rmfamily\itshape#1.}}{\qed\vspace{0.15cm}}%{\newline\vspace{0.15cm}}

\theoremstyle{thm}
\newtheorem{theorem}{Theorem}[section]
\newtheorem{lemma}[theorem]{Lemma}
\newtheorem{proposition}[theorem]{Proposition}

\newtheorem*{Main Theorem}{Main Theorem.}
\newtheorem{corollary}[theorem]{Corollary}
\newtheorem*{special theorem}{Lindeberg-Feller Theorem for Martingales}

\newtheorem{assumption}[theorem]{Assumption}

\theoremstyle{def}

\theoremstyle{rem}
\newtheorem{remark}[theorem]{{\bfseries Remark}}

\numberwithin{equation}{section}

\renewcommand{\section}{\secdef\sct\sect}
\newcommand{\sct}[2][default]{%
\refstepcounter{section}
\addcontentsline{toc}{section}{{\tocsection {}{\thesection}{\!\!\!\!#1\dotfill}}{}}
\vspace{0.7cm}
\centerline{\scshape\thesection.\ #1} \nopagebreak \vspace{0.2cm}}
\newcommand{\sect}[1]{%
\vspace{0.4cm} \centerline{\large\scshape\rmfamily #1}
\vspace{0.2cm}}

\renewcommand{\subsection}{\secdef\subsct\sbsect}
\newcommand{\subsct}[2][default]{\refstepcounter{subsection}
\addcontentsline{toc}{subsection}
{{\tocsection{\!\!}{\hspace{1.2em}\thesubsection}{\!\!\!\!#1\dotfill}}{}}
\nopagebreak\vspace{0.75\baselineskip} {\flushleft\bf
\thesubsection~\bf #1.~}
\\*[3mm]\noindent
\nopagebreak}
\newcommand{\sbsect}[1]{\nopagebreak\vspace{0.75\baselineskip}\noindent
\textbf{#1.~}\\*[3mm]\noindent
\nopagebreak}

\renewcommand{\subsubsection}{%
\secdef \subsubsect\sbsbsect}
\newcommand{\subsubsect}[2][default]{%
\refstepcounter{subsubsection} 
\addcontentsline{toc}{subsubsection}{{\tocsection{\!\!}
{\hspace{3.05em}\thesubsubsection}{\!\!\!\!#1\dotfill}}{}}
\nopagebreak
\vspace{0.15\baselineskip} \nopagebreak {\flushleft\rmfamily
\itshape\thesubsubsection
\ \rmfamily #1\/.}\ }
\newcommand{\sbsbsect}[1]{\vspace{0.1cm}\noindent
\rmfamily \itshape
\arabic{section}.\arabic{subsection}.\arabic{subsubsection} \
\sffamily #1\/.\ }

\renewcommand{\caption}[1]{%
\vglue0.5cm
\refstepcounter{figure}
\begin{minipage}{0.9\textwidth}\small {\sc Figure~\thefigure. }#1\end{minipage}}

\newcommand{\dist}{\operatorname{dist}}

\newcommand{\diam}{\operatorname{diam}}
\newcommand{\esssup}{\operatorname{esssup}}
\newcommand{\essinf}{\operatorname{essinf}}
\newcommand{\gap}{\operatorname{gap}}

\renewcommand{\d}{\operatorname{d}}

\newcommand{\textd}{\text{\rm d}\mkern0.5mu}

\newcommand{\texte}{\text{\rm e}}

\newcommand{\1}{\operatorname{\sf 1}\!}

\newcommand{\Prob}{\text{\rm Prob}}

\newcommand{\CC}{\mathcal C}

\newcommand{\HH}{\mathcal H}

\newcommand{\KK}{\mathcal K}
\newcommand{\LL}{\mathcal L}

\newcommand{\C}{\mathbb C}

\newcommand{\E}{\mathbb E}

\newcommand{\N}{\mathbb N}

\newcommand{\BbbP}{\mathbb P}

\newcommand{\R}{\mathbb R}

\newcommand{\Z}{\mathbb Z}

\newcommand{\scrF}{\mathscr{F}}

\newcommand{\scrX}{\mathscr{X}}

\newcommand{\IM}{{\Im \mathfrak m}\mkern2mu}

\newcommand{\twoeqref}[2]{(\ref{#1}--\ref{#2})}
\newcommand{\cc}{{\text{\rm c}}}

\def\myffrac#1#2 in #3{\raise 2.6pt\hbox{$#3 #1$}\mkern-1.5mu\raise 0.8pt\hbox{$#3/$}\mkern-1.1mu\lower 1.5pt\hbox{$#3 #2$}}
\newcommand{\ffrac}[2]{\mathchoice%
{\myffrac{#1}{#2} in \scriptstyle}
{\myffrac{#1}{#2} in \scriptstyle}
{\myffrac{#1}{#2} in \scriptscriptstyle}
{\myffrac{#1}{#2} in \scriptscriptstyle}
}

\newcommand{\ssup}[1] {{\scriptscriptstyle{({#1}})}}

\newcommand{\hata}{\widehat a}

\begin{document}

\vglue-0.3cm

\title[Random Schr\"odinger operators\hfill\qquad]{Eigenvalue order statistics for random Schr\"odinger\\operators with doubly-exponential tails}
\author[\qquad \hfill Biskup \and K\"onig]{Marek Biskup$^{1,2}$\, \and\,\,Wolfgang K\"onig$^{3,4}$}

\thanks{\hglue-4.5mm\fontsize{9.6}{9.6}\selectfont\copyright\,2013 by M.~Biskup and W.~K\"onig. Reproduction, by any means, of the entire article for non-commercial purposes is permitted without charge.
\vspace{2mm}
}

\maketitle

\vglue-5mm

\centerline{\textit{$^1$Department of Mathematics, UCLA, Los Angeles, California, USA}}
\centerline{\textit{$^2$School of Economics, University of South Bohemia, \v Cesk\'e Bud\v ejovice, Czech Republic}}
\centerline{\textit{$^3$Weierstra\ss-Institut f\"ur Angewandte Analysis und Stochastik, Berlin, Germany}}
\centerline{\textit{$^4$Institut f\"ur Mathematik, Technische Universit\"at Berlin, Berlin, Germany}}

\vspace{-2mm}
\begin{abstract}
We consider random Schr\"odinger operators of the form $\Delta+\xi$, where $\Delta$ is the lattice Laplacian on~$\Z^d$ and $\xi$ is an i.i.d.\ random field, and study the extreme order statistics of the eigenvalues for this operator restricted to large but finite subsets of~$\Z^d$. We show that for~$\xi$ with a doubly-exponential type of upper tail, the upper extreme order statistics of the eigenvalues falls into the Gumbel max-order class. The corresponding eigenfunctions are exponentially localized in regions where~$\xi$ takes large, and properly arranged, values. A new and self-contained argument is thus provided for Anderson localization at the spectral edge which permits a rather explicit description of the shape of the potential and the eigenfunctions. Our study serves as an input into the analysis of an associated parabolic Anderson problem.
\end{abstract}

\section{Introduction and Results}\label{sec-Intro}
%\vglue-7mm
\noindent
Random Schr\"odinger operators have been a focus of interest among mathematicians and mathematical physicists for several decades. A good representative class is lattice Schr\"odinger Hamiltonian $H_\xi$ that acts on $f\colon\Z^d\to\C$ as
\begin{equation}
\label{Ham}
(H_\xi f)(x):=\sum_{y\colon |y-x|=1}\bigl[f(y)-f(x)\bigr]+\xi(x)f(x),\qquad x\in\Z^d,
\end{equation}
with the potential $\{\xi(x)\}_{x\in\Z^d}$ sampled independently from a common law on~$\R$. The first term on the right is the lattice Laplacian so we may also write $H_\xi=\Delta+\xi$. Note that our sign conventions in~\eqref{Ham} are different from physics; hence our focus on the maximum of the spectrum.

Much is known (and unknown) about the spectral properties of~$H_\xi$. Our principal goal here is a description of the spectral \emph{extreme order statistics} for $H_\xi$ over large finite subsets of~$\Z^d$. More precisely, for a finite~$D\subset\Z^d$, let $H_{D,\xi}$ denote the operator~$H_\xi$ restricted to functions with Dirichlet boundary condition outside~$D$. This is a self-adjoint operator (a matrix) with real eigenvalues that we will label in decreasing order as
\begin{equation}
\label{E:1.2}
\lambda_D^{\ssup 1}(\xi)\ge\lambda_D^{\ssup2}(\xi)\ge\dots\ge\lambda_D^{\ssup{|D|}}(\xi).
\end{equation}
As is common in extreme-value theory, for a sequence~$D_L$ of finite subsets of~$\Z^d$ with $D_L\uparrow\Z^d$, we wish to identify sequences $a_L$ and~$b_L$ so that, as $L\to\infty$, the set of points
\begin{equation}
\label{E:eigenset}
\biggl\{\frac1{b_L}\bigl(\lambda_{D_L}^{\ssup k}(\xi)-a_L\bigr)\colon k=1,\dots,|D_L|\biggr\}
\end{equation}
tends in law to a non-degenerate point process on~$\R$.

Of course we cannot hope to do this just for any sequence of domains $D_L$ so we will content ourselves with domains that arise as scaled lattice versions,
\begin{equation}
\label{E:D_L}
D_L:=\{x\in\Z^d\colon x/L\in D\}=(LD)\cap\Z^d,
\end{equation}
of bounded and open sets~$D\subset\R^d$ with a rectifiable boundary $\partial D$. We will use $\mathfrak D$ denote the collection of all such sets. For reasons to be explained later, we will also limit ourselves to potentials whose upper tails are close to the doubly-exponential distribution,
\begin{equation}
\label{tail}
\Prob\bigl(\xi(0)>r\bigr)=\exp\bigl\{-\texte^{r/\rho}\bigr\},
%,\qquad r\to\infty,
\end{equation}
where~$\rho\in(0,\infty)$. The specific class of potentials we will consider is determined by:

\begin{assumption}
\label{ass}
Suppose $\esssup\xi(0)=\infty$ and let
\begin{equation}
F(r):=\log\log\bigl(\,\BbbP(\xi(0)>r)^{-1}\bigr),\qquad r>\essinf\xi(0).
\end{equation}
%$F(r):=\log\log\bigl(\,\BbbP(\xi(0)>r)^{-1}\bigr)$ with domain $r>\essinf\xi(0)$. 
We assume that $F$ is continuously differentiable on its domain and there is $\rho\in(0,\infty)$ such that
\begin{equation}
\label{E:4.2}
\lim_{r\to\infty}F'(r)=\frac1\rho.
\end{equation}
\end{assumption}
Our results will address not only the eigenvalues but also the associated eigenfunctions. For this, let $\{\psi_{D,\xi}^{\ssup k}\colon k=1,\dots,|D|\}$ denote an orthonormal basis of eigenfunctions of~$H_{D,\xi}$ in a (finite) set~$D\subset\Z^d$ that are labeled such that
\begin{equation}
%\label{}
H_{D,\xi}\psi_{D,\xi}^{\ssup k}=\lambda_D^{\ssup k}(\xi)\psi_{D,\xi}^{\ssup k},\qquad k=1,\dots,|D|.
\end{equation}
Although the eigenfunctions are not uniquely determined when two of the eigenvalues coincide, we can and will take these real valued.
We then define $X_k=X_k(\xi)$ by
\begin{equation}
\label{E:1.9w}
\bigl|\psi_{D,\xi}^{\ssup k}(X_k)\bigr|=\max_{x\in D}\bigl|\psi_{D,\xi}^{\ssup k}(x)\bigr|,\qquad k=1,\dots,|D|,
\end{equation}
resolving ties using the lexicographic order on~$\Z^d$. 

The family $\{(\lambda_{D_L}^{\ssup k}(\xi),\psi_{D_L,\xi}^{\ssup k})\colon k=1,\dots,|D_L|\}$ thus identifies a random sequence of points $\mathfrak a_1,\dots,\mathfrak a_{|D_L|}\in\R^{d+1}$, where $\mathfrak a_k$ is defined by
\begin{equation}
%\label{}
\mathfrak a_k:=\biggl(\,\frac {X_k(\xi)}L, \frac1{b_L}\bigl(\lambda_{D_L}^{\ssup k}(\xi)-a_L\bigr)\biggr),\qquad k=1,\dots,|D_L|.
\end{equation}
These are ordered by their last coordinate. 
This sequence induces a random measure~$\scrX_L$ on $\R^{d+1}$ by setting, for any Borel set~$B\subset\R^{d+1}$,
\begin{equation}
\label{E:1.10q}
\scrX_L(B):=\sum_{k=1}^{|D_L|}\1_{\{\mathfrak a_k\in B\}}.
\end{equation}
Note that $\scrX_L(B)$ takes values in $\N_0=\{0,1,2,\dots\}$, and $\scrX_L$ is supported on $D\times\R$.

Let us briefly recall  some facts about point processes. First, by a point process $\scrX$ on~$\R^{d+1}$ we will mean a random $\N_0\cup\{\infty\}$-valued measure on Borel sets in~$\R^{d+1}$ such that $\scrX(C)<\infty$ a.s.\ for any compact $C\subset\R^{d+1}$. The space of such measures, endowed with the topology of vague convergence, is a Polish space so convergence in law can be defined accordingly. For instance,~$\scrX_L$ converges in law to $\scrX$ if (and only if) for any continuous and compactly-supported function $f\colon\R^{d+1}\to\R$, the integral of~$f$ against $\scrX_L$ converges in law to the integral of~$f$ against~$\scrX$. 

A process $\scrX$ is a Poisson point process with intensity measure~$\mu$ if and $\scrX(B)$ is a Poisson random variable with parameter~$\mu(B)$ for any Borel set $B$ and $\scrX(B_1),\dots,\scrX(B_n)$ are independent for any pairwise disjoint Borel sets $B_1,\dots,B_n$. The principal result of the present paper is then:

\begin{theorem}[Poisson convergence; eigenfunction localization]
\label{thm-orderstat}
Fix~$d\ge1$ and let $(\xi(x))_{x\in\Z^d}$ be i.i.d.\ random variables satisfying Assumption~\ref{ass} with some $\rho\in(0,\infty)$. Then there is a sequence~$a_L$ with asymptotic growth
\begin{equation}
%\label{}
a_L=\bigl(\rho+o(1)\bigr)\log\log L,\qquad L\to\infty,
\end{equation}
such that, for any~$D\in\mathfrak D$ with scaled lattice version~$D_L$ and any choice of the normalized eigenfunctions $\{\psi_{D_L,\xi}^{\ssup k}\}$ as above, we have:
\settowidth{\leftmargini}{(1111)}
%\settowidth{\leftmargini}{(11)}
\begin{enumerate}
\item[(1)] (Eigenfunction localization) For any~$k\in\N$ and any $r_L\to\infty$,
\begin{equation}
\label{E:1.7}
\sum_{z\colon|z-X_k(\xi)|\le r_L}\bigl|\psi_{D_L,\xi}^{\ssup k}(z)\bigr|^2\,\underset{L\to\infty}\longrightarrow\,1
\qquad\text{\rm in~$\BbbP$-probability}.
\end{equation}
\item[(2)] (Poisson convergence) The process $\scrX_L$, defined via \eqref{E:1.10q} by the points
\begin{equation}
\label{E:1.8}
\biggl\{\Bigl(\frac{X_k(\xi)}L,\,\frac1\rho\bigl(\lambda_{D_L}^{\ssup k}(\xi)-a_L\bigr)\log|D_L|\Bigr)\colon k=1,\dots,|D_L|\biggr\},
\end{equation}
converges in law to the Poisson point process on $D\times\R$ with intensity measure $\textd x\otimes\texte^{-\lambda}\textd\lambda$.
\end{enumerate}
\end{theorem}

Restricting the points in $\scrX_L$ to the last coordinate shows that the set \eqref{E:eigenset} converges in law to a Poisson process provided $a_L$ is as above and $b_L\asymp\log L$. The limit law can be described concisely and explicitly as follows:

\begin{corollary}[Eigenvalue order-statistics]
\label{cor-1.3}
Assume the setting of Theorem~\ref{thm-orderstat} and let $D_L$ be as in \eqref{E:D_L}. Then the (upper) order statistics of the eigenvalues lies in the domain of attraction of the Gumbel universality class. In particular, for each~$k\in\N$, 
\begin{multline}
%\label{}
\qquad\quad
\label{eigenvalueorderstat}
\bigl(\texte^{-\frac1\rho(\lambda_{D_L}^{\ssup 1}(\xi)-a_L)\log|D_L|},\dots,\texte^{-\frac1\rho(\lambda_{D_L}^{\ssup k}(\xi)-a_L)\log|D_L|}
\bigr)\,
\\
\overset{\text{\rm law}}{\underset{L\to\infty}\longrightarrow}\,\,\,\bigl(Z_1,Z_1+Z_2,\dots,Z_1+\dots+Z_k),
\qquad\quad
\end{multline}
where $Z_1,Z_2,\dots$ are i.i.d.\ exponential with parameter~one. Equivalently, the vector on the left tends in law to the first $k$ points of a Poisson point process on $[0,\infty)$ with intensity one.
\end{corollary}

We remind the reader that the Gumbel universality class is one out of three possible non-degenerate limit distributions for order statistics of i.i.d.\ random variables; see e.g.\ de Haan and Ferreira~\cite{deHaan}. Gumbel is also the extreme-order class associated with the doubly exponential tails (see Section~\ref{sec:extreme-val}). However, as many values of the field need to ``cooperate'' to create conditions for an extremal eigenvalue, the shift~$a_L$ required for the~$\xi$'s differs from the one above by a non-vanishing amount. Explicitly, under Assumption~\ref{ass} and for any  $D\in\mathfrak D$,
\begin{equation}
\label{E:1.16}
\max_{x\in D_L}\xi(x)-\lambda_{D_L}^{\ssup 1}(\xi)\,\underset{L\to\infty}\longrightarrow\,\chi
\qquad\text{\rm in~$\BbbP$-probability},
\end{equation}
where $\chi=\chi(\rho,d)$ is the quantity in $(0,2d]$ given by
\begin{equation}
\label{chi-def}
\chi:=-\sup\bigl\{\lambda^{\ssup 1}(\varphi)\colon\varphi\in\R^{\Z^d},\,\LL(\varphi)\le1\bigr\},
\end{equation}
while $\LL\colon\R^{\Z^d}\to(0,\infty]$ is defined as
\begin{equation}
\label{LL-phi}
\LL(\varphi):=\sum_{x\in\Z^d}\texte^{\varphi(x)/\rho}.
\end{equation}
We note that $\lambda^{\ssup 1}(\varphi)$, the supremum of the spectrum of $\Delta+\varphi$, is an isolated (simple) eigenvalue for all potentials~$\varphi$ with $\LL(\varphi)<\infty$.

The limit \eqref{E:1.16} is actually well known (albeit under different assumptions) from earlier studies of $H_\xi$ with $\xi$ having doubly-exponential tails; see Section~\ref{sec-connections}. The function $\varphi\mapsto\LL(\varphi)$ is encountered in these studies as well; it is the large-deviation rate function for the field~$\xi$ and so it plays an important role in estimating the probability that~$\xi$ exceeds a given function in a given domain; see~\eqref{E:4.11} for an explicit statement in this vain.

The proof of the limiting Poisson statistics \eqref{eigenvalueorderstat} is based on constructing a coupling to i.i.d.\ random variables; cf Theorem~\ref{prop-4.6} for a precise formulation. This coupling applies to a whole non-degenerate interval at the top of the spectrum. The statement \eqref{E:1.7} implies exponential localization of leading eigenfunctions at the lattice scale. We state a quantitative decay bound that concerns the leading eigenfunctions:

\begin{theorem}
\label{thm2.4}
For $X_k(\xi)$ as in Theorem~\ref{thm-orderstat} and each $k\ge1$ the following holds with probability tending to one as $L\to\infty$: There exist (deterministic) constants $c_1,c_2>0$ such that
\begin{equation}
%\label{}
\bigl|\psi_{D_L,\xi}^{\ssup k}(z)\bigr|\le c_1\texte^{-c_2|z-X_k(\xi)|},\qquad z\in D_L.
\end{equation}
Moreover, for larger separations from $X_k(\xi)$ we in fact get
\begin{equation}
\label{E:1.20q}
\bigl|\psi_{D_L,\xi}^{\ssup k}(z)\bigr|\le c_1'\texte^{-c_2'(\log\log L)|z-X_k(\xi)|},\qquad \bigl|z-X_k(\xi)\bigr|\ge\log L.
\end{equation}
for some non-random constants $c_1',c_2'>0$.
\end{theorem}

We emphasize that the methods of the present paper are largely independent of the existing techniques for proving Anderson localization. In particular, our approach permits a rather explicit characterization of the location, size and shape of the potential and the corresponding eigenfunction at (and near) the top of the spectrum. 

\section{Road map to proofs}
\label{sec2}\noindent
We proceed to discuss the key ideas of the proofs. We break the main argument into a sequence of stand-alone steps which, we believe, are of independent interest. 

As already alluded to, our results are a manifestation of Anderson localization, discovered in 1958 by Anderson~\cite{Anderson} and studied extensively by mathematicians in the past four decades. The word ``localization'' refers to the fact that, when the fields~$\xi$ are non-degenerate i.i.d.\ random variables, the spectrum of the Hamiltonian $H_\xi$ in \eqref{Ham} will contain a band of proper eigenvalues with exponentially {localized} eigenfunctions. (This is in contrast with the situation when~$\xi$ is only periodic, where the spectrum has a band structure but remains continuous, by the classic Bloch theory.) We refer to, e.g., Pastur and Figotin~\cite{Pastur-Figotin}, Stollmann~\cite{Stollmann}, Carmona and Lacroix~\cite{Carmona-Lacroix} and Hundertmark~\cite{Hundertmark} for further details and explanations.

Invariably, all existing proofs of Anderson localization are based on controlling the Green function associated with~$H_{D,\xi}$,
\begin{equation}
%\label{}
G_{D,\xi}(x,y;z):=\bigl\langle\delta_x, (z-H_{D,\xi})^{-1}\delta_y\bigr\rangle,\qquad \IM\, z>0,
\end{equation}
in the limit as~$z$ ``radially'' approaches the real line from the upper half plane in~$\C$. The aim is to show that $G_{D,\xi}(x,y;z)$ exhibits exponential decay in~$|x-y|$ in the said limit; a key challenge is to avoid~$z$ hitting an eigenvalue where the constant in front of the exponentially decaying term becomes divergent. 
Various averaging methods have been developed with this purpose in mind. For instance, the fractional-moment method of Aizenman and Molchanov~\cite{Aizenman-Molchanov} generally yields bounds of the form
\begin{equation}
\label{E:FFB}
\E\bigl(|G_{D,\xi}(x,y;z)|^s\bigr)\le c_1\texte^{-c_2|x-y|}
\end{equation}
uniformly in $D$ and $\IM(z)>0$, where $c_1,c_2\in(0,\infty)$ are constants while the exponent $s\in(0,1)$ is tied to H\"older continuity of the probability density of~$\xi(0)$.  Argument from spectral theory for \emph{infinite-volume} operators (perfected into the so called Simon-Wolff criteria~\cite{Simon-Wolff}) then permit one to infer from \eqref{E:FFB} the existence of eigenvalues with exponentially decaying eigenfunctions. Approaches based on ``finite-volume criteria'' also exist (Aizenman, Schenker, Friedrich and Hundertmark~\cite{ASFH}), but they are still versed in the language of the Green function.

Our approach is different from the above in a number of important aspects. We work directly with individual eigenvalues in a finite volume and control their dependence on the configuration of random fields. This permits us to characterize geometrically the regions where the eigenfunctions are localized. On the technical side, we manage to avoid working with complex weights and the Green function. Large-deviation theory naturally lurks in the background although, for the most part, we proceed by direct estimates. Although our method is, in its present form, tailored to the study of the upper edge of the spectrum for operators $H_{D,\xi}$ in finite~$D$ with unbounded i.i.d.\ random fields, we believe that an extension for bounded fields is possible. 

We will now proceed to describe the main steps of our approach formulating the key parts thereof as separate theorems.

\subsection*{STEP 1: Domain truncation and component trimming}
A good deal of our proof of Theorem~\ref{thm-orderstat} focuses on individual eigenvalues.
The starting observation is that the field configuration $\xi$ in regions where $\xi$ is smaller than an eigenvalue $\lambda$ is of little relevance for~$\lambda$. For~$A>0$ and $R\in\N$, consider the set
\begin{equation}
\label{E:2.3}
D_{R,A}(\xi):=\bigcup_{z\in D\colon\xi(z)\ge\lambda_D^{\ssup 1}(\xi)-2A}
B_R(z)\cap D,
\end{equation}
where $B_R(z):=\{x\in\Z^d\colon |z-x|_1\le R\}$. We will occasionally refer to the field values in $D_{R,A}$ ``large'' while those not in this set as ``small.'' Deterministic arguments show:

\begin{theorem}[Domain truncation]
\label{thm-truncation}
Let~$A>0$ and $R\in \N$ be such that $2d(1+\frac A{2d})^{1-2R}\le\frac A2$. Then for any $\xi$, any $U$ with $D_{R,A}(\xi)\subset U\subset D$ and any~$k\in\{1,2,\dots,|U|\}$ such that
\begin{equation}
\label{lambda-theta-xi-initial}
\lambda_D^{\ssup k}(\xi)\ge\lambda_D^{\ssup 1}(\xi)-\frac A2
\end{equation}
we have
\begin{equation}
\label{eigen-distance}
\bigl|\lambda_D^{\ssup k}(\xi)-\lambda_U^{\ssup k}(\xi)\bigr|\le 2d\Bigl(1+\frac A{2d}\Bigr)^{1-2R}.
\end{equation}
\end{theorem}

The vehicle that brings us to this conclusion is, not surprisingly, analysis of rank-one perturbations of~$H_{D,\xi}$. However, unlike for the corresponding arguments in the proofs of, e.g., \eqref{E:FFB}, where the perturbation occurs at a single vertex, here we address single eigenvalues (rather than the Green function) and we perturb the configuration in all of~$D\setminus U$. 

In the specific context of doubly-exponential tails, we will use the conclusion in \eqref{eigen-distance} for~$D$ replaced by $D_L$, the cutoff~$A$ fixed to a small number (less than~$\chi$) and $R$ growing slowly to infinity with~$L$. Under such conditions, the components of~$U:=D_{R,A}$ become very sparse and their geometry can be analyzed by straightforward estimates. In particular, due to Dirichlet boundary conditions, the spectrum of $H_{U,\xi}$ is the union of the spectra in the connected components of~$U$. 

If~$C$ is such a component and $\lambda_C^{\ssup 1}(\xi)<\lambda_D^{\ssup1}(\xi)-A$, then~$C$ cannot contribute to the set of eigenvalues covered by \eqref{eigen-distance}. We can thus remove $C$ from~$U$ and still maintain the control provided by Theorem~\ref{thm-truncation}. This permits systematic component ``trimming'' that helps significantly reduce the number of connected components of concern.

\subsection*{STEP 2: Reduction to one eigenvalue per component}
The removal of irrelevant components of $D_{R,A}$ eliminates much of the geometric complexity of the underlying field configuration. An issue that comes up next is what part of the spectrum in each component needs to be taken into account. Here we will observe that in components that have a chance to contribute, all but the leading eigenvalue can safely be disregarded. This effectively couples the top of the spectrum in~$D$ to the set of principal eigenvalues in the components of~$D_{R,A}$.

For~$C\subset\Z^d$ finite, consider the finite-volume version of~\eqref{LL-phi},
\begin{equation}
%\label{}
\LL_C(\varphi):=\sum_{x\in C}\texte^{\varphi(x)/\rho}.
\end{equation}
Similarly, we introduce the finite-volume analogue of \eqref{chi-def}, 
\begin{equation}
\label{chi-C-alter}
\chi_C:=-\sup\bigl\{\lambda_C^{\ssup 1}(\xi)\colon \xi\in\R^{\Z^d},\,\LL_C(\xi)\le1\bigr\}.
\end{equation}
Then we have the following deterministic estimate:

\begin{proposition}[Spectral gap]
\label{prop-eigengap}
Let~$C\subset\Z^d$ be finite. If for some $K\ge0$,
\begin{equation}
\label{eigen-gap}
\lambda_C^{\ssup 1}(\xi)-\lambda_C^{\ssup 2}(\xi)\le K,
\end{equation}
then also
\begin{equation}
\label{eigen-chi-gap}
\lambda_C^{\ssup 1}(\xi)-\rho\log\LL_C(\xi)\le-\chi_C+K-\rho\log2.
\end{equation}
\end{proposition}
 
Our use of this proposition requires observations from large-deviation theory for double-exponential i.i.d.\ random fields: First, for potentials satisfying Assumption~\ref{ass}, $\LL_C$ acts as a large-deviation rate function for finding a specific potential profile in~$C$. More explicitly, in Section~\ref{sec:extreme-val} we will show that, for some $\hata_L\to\infty$,
\begin{equation}
\label{5.1-eq}
\BbbP(\xi(0)\ge\hata_L+s)=\frac1{L^{d\theta(1+o(1))}}\quad\text{where}\quad\theta:=\texte^{s/\rho},
\end{equation}
with $o(1)\to0$ as $L\to\infty$ uniformly in compact sets of~$s$. Thus, for any given $\varphi\colon C\to\R$,
\begin{equation}
\label{E:4.11}
\BbbP\bigl(\xi\ge\hata_L+\varphi\text{ in }C\bigr)=L^{-d\LL_C(\varphi)[1+o(1)]},\qquad L\to\infty.
\end{equation}
The set inclusion
\begin{equation}
\label{E:recall-union2}
\bigl\{\xi\colon\lambda_C^{\ssup 1}(\xi)\ge a\bigr\}\subseteq\bigcup_{\varphi\colon\LL_C(\varphi)\ge1}\bigl\{\xi\colon\xi\ge a+\chi_C+\varphi\text{ \rm on }C\bigr\},
\end{equation}
which is derived readily from the alternative formula for $\chi_C$, 
\begin{equation}
\label{chi-C-def}
\chi_C=-\sup_\xi\bigl[\lambda_C^{\ssup 1}(\xi)-\rho\log\LL_C(\xi)\bigr],
\end{equation}
and a union bound show that large eigenvalues in $D_L$ will thus come only with potential profiles for which $\varphi:=\xi-\hata_L-\chi_C$ obeys $\LL_C(\varphi)\approx1$.

Returning to the role of Proposition~\ref{prop-eigengap} in our proofs, we note that its main conclusion can be used to derive
\begin{equation}
\label{3:23-eq}
\begin{aligned}
\lambda_C^{\ssup 1}(\xi)\ge a'&\,\,\text{ \rm\small AND }\,\,\lambda_C^{\ssup 1}(\xi)-\lambda_C^{\ssup 2}(\xi)\le\tfrac12\rho\log2\qquad\\
&\Longrightarrow\qquad\xi\ge\varphi+a+\chi_C\,\text{ \rm in $C$ for some $\varphi$ satisfying }\,\LL_C(\varphi)\ge u,
\end{aligned}
\end{equation}
where $u$ is defined by
\begin{equation}
\label{3:23-eq2}
\log u:=\frac{a'-a}\rho+\frac12\log2.
\end{equation}
Since $u>0$ for~$a'> \hata_L-\chi-\frac\rho2\log2$, the event that there is~$C$ satisfying the conditions on the left of \eqref{3:23-eq} has probability $o(L^{-d})$ and thus will not occur once~$L$ is sufficiently large. It follows that, whenever~$\lambda_C^{\ssup 1}(\xi)$ for some component~$C$ in~$D_L$ is close to its optimal value,~$\lambda_C^{\ssup 2}(\xi)$ is at least a deterministic constant below~$\lambda_C^{\ssup 1}(\xi)$. In short, only the top eigenvalue in each connected component of $D_{R,A}$ need be considered.

\subsection*{STEP 3: Coupling to i.i.d.\ variables}
Our previous observations permit us to design a coupling between the eigenvalues in a small interval near the top the spectrum of $H_{D,\xi}$ and a family of i.i.d.\ random variables. There is a number of ways how such a coupling can be formulated; we will present one that is based on a regular partition of~$\Z^d$ into square boxes. 

As we will work, from now on, with~$L$-dependent objects, we need to fix two sequences that determine the main scales of the problem: a sequence of integers~$R_L$ satisfying
\begin{equation}
\label{E:RLcond}
\frac{R_L}{\log\log L}\,\underset{L\to\infty}\longrightarrow\,\infty\quad\text{but}\quad R_L=(\log L)^{o(1)},
\end{equation}
which will govern the size of the connected components and spatial range of the perturbation arguments described above, and a sequence of integers~$N_L$ such that 
\begin{equation}
\label{E:RLNL}
%N_L=(\log L)^{o(1)}\quad\text{and}\quad 
\frac{N_L}{R_L}\,\underset{L\to\infty}\longrightarrow\,\infty\qquad\text{and}\qquad \limsup_{L\to\infty}\,\frac{\log N_L}{\log L}< 1 %\frac{d}{d+1},
\end{equation}
which determines the size of the boxes in the partition.

Using a natural partition of~$\Z^d$ into square boxes of side~$N_L+1$, for each such box $B_{N_L+1}$, consider the sub-box of side~$N_L$ induced by the embedding $B_{N_L}\subset B_{N_L+1}$. We will call these $N_L$-boxes. For $D\in\mathfrak D$, let $D_L$ be its scaled lattice version \eqref{E:D_L} and let $B^{(i)}_{N_L}$, $i=1,\dots,m_L$, denote the collection of those $N_L$-boxes that are entirely contained in~$D_L$. Since $\partial D$ is rectifiable, we have
\begin{equation}
%\label{}
|D_L|=m_L\,N_L^d\bigl(|D|+o(1)\bigr),
\end{equation}
where $|D|$ denotes the Lebesgue volume of~$D$. By \eqref{E:RLNL}, $m_L\to\infty$ as $L\to\infty$.

Given a configuration~$\xi$, let $\lambda_i(\xi)$ be a shorthand for the principal eigenvalue in $B^{(i)}_{N_L}$,
\begin{equation}
\label{E:4.27c}
\lambda_i(\xi):=\lambda_{B^{(i)}_{N_L}}^{\ssup 1}(\xi),\qquad i=1,\dots,m_L.
\end{equation}
Since these eigenvalues depend on disjoint subsets of the~$\xi$ variables, the random variables $\{\lambda_i(\xi)\colon i=1,\dots,m_L\}$ are i.i.d.

\begin{theorem}[Coupling to i.i.d.\ process]
\label{prop-4.6}
Let $D\in\mathfrak D$.
For sequences $(R_L)$, $(N_L)$ and $D_L$ as above, let $\hat\lambda_1(\xi),\dots,\hat\lambda_{m_L}(\xi)$ be the sequence $\lambda_1(\xi),\dots,\lambda_{m_L}(\xi)$ from \eqref{E:4.27c} listed in decreasing order. Under Assumption~\ref{ass}, there is an $A>0$ such that the event
\begin{equation}
\label{E:2.20ww}
\forall k\in\{1,\dots,m_L\}\colon\quad\hat\lambda_1(\xi)-\hat\lambda_k(\xi)<A\,\,\,\Rightarrow\,\,\,
\bigl|\lambda_{D_L}^{\ssup k}(\xi)-\hat\lambda_k(\xi)\bigr|<4d\Bigl(1+\frac A{2d}\Bigr)^{1-2R_L}
\end{equation}
occurs with probability tending to one as $L\to\infty$.
\end{theorem}

The proof is based on Theorem~\ref{thm-truncation} and estimates of the probability that a component with an appreciable principal eigenvalue intersects the boundary of a partition box, or lies within~$N_L$ of the boundary of~$D_L$. (These are the reason for the restrictions in \eqref{E:RLNL}.)

\subsection*{STEP 4: Identifying max-order class}
\noindent
Theorem~\ref{prop-4.6} brings the proof of the Poisson statistics in Theorem~\ref{thm-orderstat} and Corollary~\ref{cor-1.3} to the realm of standard extreme-order limit theory (see e.g.\ de Haan and Ferreira~\cite{deHaan}). Naturally, one starts by defining the centering sequence $a_L$ as
\begin{equation}
\label{E:4.13q}
\BbbP\bigl(\lambda_{B_{N_L}}^{\ssup 1}\ge a_L\bigr)=\Bigl(\frac{N_L}L\Bigr)^d.
\end{equation}
(Such an $a_L$ exists because the $\xi$'s, and also the $\lambda_{B_{N_L}}^{\ssup 1}(\xi)$'s, are continuously distributed.)
In order to determine the correct scaling of the spacings between the eigenvalues, and generally place the limit distribution in the Gumbel max-order class, it remains to prove:

\begin{theorem}[Max-order class of local eigenvalues]
\label{thm-maxorder}
Suppose Assumption~\ref{ass} and let~$b_L$ obey
\begin{equation}
\label{E:4.46q}
b_L\log L\,\underset{\,L\to\infty\,}\longrightarrow\,\frac \rho d.
\end{equation}
Then for each~$s\in\R$, 
\begin{equation}
\label{4.16-eq}
\BbbP\bigl(\lambda_{B_{N_L}}^{\ssup 1}\ge a_L+sb_L\bigr)=\texte^{-s}\,\Bigl(\frac{N_L}{L}\Bigr)^d\bigl(1+o(1)\bigr),
\end{equation}
with $o(1)\to0$ as~$L\to\infty$ uniformly on compact sets of~$s$.
\end{theorem}

The proof is based on regularity of the probability density of~$\xi$ supplied by Assumption~\ref{ass}. In simple terms, the change of the eigenvalue by $sb_L$ can be achieved by shifting the whole $\xi$ configuration by the same amount. A  catch is that this would be too costly (i.e., inefficient) to perform in the entire $N_L$-box; rather one has to do this only in those parts of the box where the relevant contribution comes from. 

Before we proceed to the next step, note that, on the basis of Theorems~\ref{prop-4.6} and~\ref{thm-maxorder}, we are already able to conclude the limit statement~\eqref{eigenvalueorderstat}. (The condition on the left of \eqref{E:RLcond} ensures that the error in \eqref{E:2.20ww} is much smaller than the spacing between eigenvalues.)

\subsection*{STEP 5: Eigenfunction localization}
The final task before us is a control of the spatial localization of the eigenfunctions. As is well known, the main obstruction to localization is degeneracy of eigenvalues. Two techniques exist for dealing with this problem: averaging and multiscale analysis. In our context, we are able to address the problem directly by developing a deterministic link between the spatial decay of an eigenfunction and the distance of the associated eigenvalue to other eigenvalues. A key novel fact is that the (still needed) non-degeneracy of the eigenvalues will be supplied by the already-proved extreme-order limit theorem.

Let $\mathfrak C_{R,A}:=\mathfrak C_{R,A}(\xi)$ denote the set of connected components of $D_{R,A}(\xi)$ and, for $V\subset\Z^d$, let~$\partial V$ mark the set of vertices outside $V$ that have an edge into~$V$. 
We will measure the decay of the eigenfunctions in terms of a distance-like object $\d(x,\CC)$, indexed by vertices~$x\in\Z^d$ and components~$\CC\in\mathfrak C_{R,A}$, on which we impose the following requirements:
\settowidth{\leftmargini}{(111111)}
\begin{enumerate}
\item[(D0)]
$\d(z,\CC)\ge0$ for all $z$ and~$\CC$ with $\d(z,\CC)=0$ whenever $z\in\CC$.
\item[(D1)]
For all $z\in D\setminus D_{R,A}(\xi)$, all $ y\in\partial B_R(z)$ and all $\CC\in\mathfrak C_{R,A}$, we have $\d(z,\CC)\le\d(y,\CC)+R$.
\item[(D2)]
For all $\CC'\ne\CC$, all $ z\in\CC'$ and all $ y\in\partial\CC'$ we have $\d(z,\CC)\le \d(y,\CC)+1$.
\end{enumerate}
An example of such $\d(\cdot,\cdot)$ is constructed as follows: Define a graph by contracting all components in $\mathfrak C_{R,A}$ to a single vertex while keeping the edges between the (new) vertex corresponding to component $\CC$ and all (old) vertices on $\partial\CC$ --- which, by the fact that $\CC$ is in~$\mathfrak C_{R,A}$ do not lie in another component in~$\mathfrak C_{R,A}$. Then set $\d(z,\CC)$ to the corresponding graph-theoretical distance from~$z$ to the (vertex corresponding to) component~$\CC$.

\begin{theorem}[Eigenfunction decay]
\label{prop-exp-decay}
Assume~$R\ge1$ and $A>0$ obey $\varepsilon_R:=2d(1+\frac A{2d})^{1-2R}<A/2$.
Let~$\lambda$, resp.,~$\psi$ be a Dirichlet eigenvalue, resp., a corresponding eigenfunction of~$H_{D,\xi}$ such that 
\begin{equation}
\label{E:3.43}
\lambda\ge\lambda_D^{(1)}(\xi)-\frac A2+\epsilon_R.
\end{equation}
Assume the following:
\settowidth{\leftmargini}{(1111)}
\begin{enumerate}
\item[(1)] $\gap(\lambda)$, the distance of $\lambda$ to the nearest eigenvalue of $H_{D,\xi}$, obeys $\gap(\lambda)>10\varepsilon_R$,
\item[(2)] there is $h>0$ such that, for any self-avoiding (nearest-neighbor) path $x_1,\dots,x_R$ in~$D$,
\begin{equation}
\label{E:4.40}
\xi(x_j)<\lambda\text{ for all }j=1,\dots, R
\quad\Rightarrow\quad\prod_{j=1}^R\frac{2d}{2d+\lambda-\xi(x_j)}\le\texte^{-h R},
\end{equation}
\item[(3)]
for some $\delta\in(0,1)$, 
\begin{equation}
\label{E:4.41}
\frac{\gap(\lambda)-2\varepsilon_R}{8d}\wedge1>4\texte^{-(1-\delta)h R+\delta h}\sqrt{|\partial\CC'|},
\qquad \CC'\in\mathfrak C_{R,A}.
\end{equation}
\end{enumerate}
Then there is $\CC\in\mathfrak C_{R,A}$ such that
\begin{equation}
\label{eigenvdecay}
\bigl|\psi(z)\bigr|\le\texte^{-\delta h\d(z,\CC)},\qquad z\in D.
\end{equation}
\end{theorem}

In order to appreciate this general result, we again place ourselves in an $L$-dependent setting of domains~$D:=D_L$, with $R:=R_L$ satisfying \eqref{E:RLcond}. Under Assumption~\ref{ass} on the upper tail of~$\xi(0)$, we then have \eqref{5.1-eq}, which readily yields \eqref{E:4.40} for any given (fixed) $h>0$. Since the component sizes are at most polylogarithmic in~$R_L$, condition (1) and (3) is satisfied as soon as $\gap(\lambda)$ is larger than exponentially-small in~$R_L$. For the leading eigenvalues, the gap is at most order $(\log L)^{-1}$; thanks to \eqref{E:RLcond}, the bound  \eqref{eigenvdecay} thus applies to $\psi:=\psi_{D_L,\xi}^{\ssup k}$ for all~$k\ge1$. A straightforward comparison between $\d(z,\CC)$ and the Euclidean distance (see~Lemma~\ref{lemma-7.2}) then prove the decay estimate \eqref{E:1.20q}.

Most of what is left in this paper consists of proofs of the above claims in full technical detail. In particular, Section~\ref{sec-truncation} deals with various deterministic spectral estimates leading to the execution of Steps~1 and~2 in the above scheme. Deterministic bounds underpinning eigenfunction localization (Step~5 above) appear in Section~\ref{sec-EigenvLoc}. Coupling to i.i.d.\ random variables (Step~3) is performed in Section~\ref{sec-orderstat}. Step~4 is the subject of Section~\ref{sec6}, where we also conclude the proofs of eigenvalue order statistics. In Section~\ref{sec7}, this feeds into the proof of eigenfunction localization (the probabilistic part of Step~5) and concludes the proof of Theorems~\ref{thm-orderstat}(2) and~\ref{thm2.4}. 

\section{Connections}
\label{sec-connections}\noindent
Before we move to actual proofs, let us pause shortly to make the requisite connections to the existing literature. References have insofar been largely suppressed in other to keep the flow of explanations of results and ideas of proofs.

Attempts to describe the statistics of the spectrum of random Schr\"odinger operators are as old as the subject itself. In the localization regime, the statistics was expected to be Poisson-like. This was proved by Killip and Nakano \cite{KN07}, based on techniques developed in Aizenman and Molchanov~\cite{Aizenman-Molchanov}, Wegner~\cite{Wegner}, Minami~\cite{Minami}. In particular, they showed that, for energies where the fractional moment bound on the resolvent applies, the point process of \emph{unfolded} eigenvalues scales to a homogeneous Poisson process. Here an unfolded eigenvalue is the quantity $I(\lambda^{\ssup k}_D)$ where $I$ is the {integrated density of states} (cf~Carmona and Lacroix \cite{Carmona-Lacroix} or Veseli\'c~\cite{Veselic}). 

Killip and Nakano's result has been further extended by Germinet and Klopp with statements that apply both in the bulk of the spectrum~\cite{GK11a} and, in $d=1$, also close to the spectral edge~\cite{GK11b} (extensions to arbitrary~$d$ require a modified kinetic term). However, these results are formulated relative to a \emph{fixed} reference point in the unfolded spectrum, and they do not seem to apply in our situations where the maximal eigenvalues tend to infinity with~$L$.

Our work is closer in spirit to the studies of the potential tails that are heavier than doubly exponential; e.g., Gren\-kova, Molchanov and Sudarev~\cite{GMS1,GMS2},  Astrauskas~\cite{A1,A2}, van der Hofstad, M\"orters and Sido\-rova~\cite{HMS08} and K\"onig, Lacoin, M\"orters and Sidorova~\cite{KMS07}. Austrauskas' study~\cite{A2} includes also doubly-exponential tails \eqref{tail} --- despite his vigorous insistence on the contrary throughout the abstract and introduction --- but only for $\rho$ very large. However, in all these works the corresponding eigenfunctions are localized  more or less at a \emph{single} lattice site; namely, a high excess value of the random potential. For the doubly exponential tails with general value of~$\rho$ this is no longer the case and this is exactly what makes these tails a challenge.

We note that in G\"artner, K\"onig and Molchanov \cite{GKM07}, the asymptotic shape of the potential in the localization regions in a large box $D$ was identified as the one of the maximizers in \eqref{chi-def}. Furthermore, an explicit form of exponential localization was proved for the principal eigenfunctions of $H_{D,\xi}$, after removing the top values of $\xi$ in all the other localization regions; the shape of these eigenfunctions was identified as well. However, the method there was based on a tedious random walk enumeration technique, which we do not follow here. To keep the present paper self-contained, we also do not use any partial result of~\cite{GKM07}.

The class of doubly exponential tails was identified rather early in the studies of the parabolic Anderson problem (G\"artner and Molchanov~\cite{GM90,GM98}, G\"artner and den Hollander~\cite{GH99}, G\"artner and K\"onig~\cite{GK05} and G\"artner, K\"onig and Molchanov~\cite{GKM07}). Other classes of potential upper tails have been identified later (Biskup and K\"onig~\cite{BK01,BK01b}, van der Hofstad, K\"onig and M\"orters~\cite{HKM06}). In these cases the support of the leading eigenvalue in a set of side~$L$ grows to infinity with~$L$. As already said, although the present paper deals only with doubly-exponential tails, we believe that the bulk of the method developed in this work applies to the other cases as well.

\section{Domain truncation and spectral gap}
\label{sec-truncation}\noindent 
The goal of this section is to establish Theorem~\ref{thm-truncation} that underlies all subsequent derivations in this paper. In addition, we will prove bounds on the distance between the first and second leading eigenvalue, i.e., the spectral gap, as stated in Proposition~\ref{prop-eigengap}. 

\subsection{Martingale argument}
The proof of Theorem~\ref{thm-truncation} is based on the fact that eigenfunctions decay rapidly away from $D_{R,A}(\xi)$. We will control the rate of this decay by a martingale argument. Let $Y:=(Y_k)_{k\in\N_0}$  denote a discrete-time simple symmetric random walk on $\Z^d$. We will write~$P^x$, resp., $E^x$ for the law, resp., expectation for the walks starting from $x\in\Z^d$ and let $\scrF_n:=\sigma(Y_0,\dots,Y_n)$ denote the canonical filtration associated with~$Y$.

\begin{lemma}
\label{lemma-localization}
Let~$\lambda:=\lambda_D(\xi)$, resp., $\psi:=\psi_{D,\xi}$ be a Dirichlet eigenvalue, resp., a corresponding eigenfunction of $H_{D,\xi}$. Define
\begin{equation}
%\label{}
\tau:=\inf\bigl\{k\in\N_0\colon \xi(Y_k)\ge\lambda\text{ \rm\small OR }Y_k\not\in D\bigr\}
\end{equation}
and denote $M_0:=\psi(Y_0)$ and, for $1\le n\le\tau$,
\begin{equation}
\label{Mn-def}
M_n:=\psi(Y_n)\,\prod_{k=0}^{n-1}\,\frac{2d}{2d+\lambda-\xi(Y_k)}.
\end{equation}
Then, under $P^x$ for any $x\in\Z^d$, the process $M^\tau=(M_{\tau\wedge n})_{n\in\N_0}$ is a martingale for the canonical filtration $(\scrF_n)_{n\in\N_0}$.
\end{lemma}

\begin{proofsect}{Proof} If $x\not\in D$, then $\tau=0$ and $M_{\tau\wedge n}=0$ for any $n$, $P^x$-a.s. For $x\in D$ the following holds $P^x$-a.s.: On~$\{\tau\ge n\}$ we have $\lambda-\xi(Y_k)\ge 0$ for $k\le n-1$ and hence $|M_n|\leq |\psi(Y_n)|\leq \max_x|\psi(x)|$; i.e.,~$M_n$ is bounded. On~$\{\tau>n\}$, the conditional expectation of~$\psi(Y_{n+1})$ given~$Y_0,\dots,Y_n$ equals $\psi(Y_n)+\frac1{2d}(\Delta \psi)(Y_n)$. Writing this using $(\Delta+\xi)\psi=\lambda \psi$ shows that, on~$\{\tau>n\}$,
\begin{equation}
\begin{aligned}
E^x(M_{n+1}|Y_0,\dots,Y_n)&=E^x(\psi(Y_{n+1})|Y_0,\dots,Y_n)\prod_{k=0}^{n}\frac{2d}{2d+\lambda-\xi(Y_k)}
\\
&=\Bigl[\psi(Y_n)+\frac1{2d}(\Delta \psi)(Y_n)\Bigr]\prod_{k=0}^{n}\frac{2d}{2d+\lambda-\xi(Y_k)}
\\
&=\psi(Y_n)\Bigl[1+\frac1{2d}\bigl(\lambda-\xi(Y_n)\bigr)\Bigr]\prod_{k=0}^{n}\frac{2d}{2d+\lambda-\xi(Y_k)}=M_n.
\end{aligned}
\end{equation}
It follows that $M^\tau$ is a martingale.
\end{proofsect}

The next lemma expresses the desired consequence of the martingale property. (The set $D'$ will later be taken to be $D\setminus U$ with $U$ as in Theorem~\ref{thm-truncation}.)

\begin{lemma}
\label{lemma-2-norm}
Let $\lambda$ and~$\psi$ be as in Lemma~\ref{lemma-localization}. Given~$A'\ge A>0$ and $R\in\N$, let~$D'\subset D$ be such that $\xi\leq\lambda-A' $ on $D'$ and 
\begin{equation}
\label{D-prime-n}
x\in D\quad\text{ \small\rm AND }\quad\xi(x)\ge\lambda-A\qquad\Longrightarrow\qquad\dist(x,D')\ge R,
\end{equation}
where ``{\rm dist}'' denotes the $\ell^1$-distance on~$\Z^d$. Then
\begin{equation}
\label{E:3.8}
\sum_{x\in D'}\bigl|\psi(x)\bigr|^2\le \Bigl(1+\frac A{2d}\Bigr)^{2-2R}\Bigl(1+\frac{A'}{2d}\Bigr)^{-2}\Vert \psi\Vert_2^2.
\end{equation}
\end{lemma}

\begin{proofsect}{Proof}
As the square of a bounded martingale is a submartingale, we have
\begin{equation}
\label{submartingale}
\bigl|\psi(x)\bigr|^2=E^x |M_{\tau\wedge 0}|^2\le E^x |M_{\tau\wedge R}|^2,\qquad x\in D'.
\end{equation}
By our assumptions on~$D'$, any path of the simple random walk started at~$x\in D'$ will either leave~$D$ or stay in the region where~$\xi<\lambda-A$ for at least~$R-1$ steps. This implies that, on the event $\{\tau<R\}$, we  necessarily have $Y_{\tau\wedge R}=Y_\tau\not\in D$ with $P^x$-probability one. Hence, $M_{\tau\wedge R}=0$ on $\{\tau<R\}$ and so $E^x(\1\{\tau<R\} |M_{\tau\wedge R}|^2)=0$.

On the other hand, on $\{\tau\geq R\}$, each term in the product in \eqref{Mn-def} is bounded by~$(1+\frac A{2d})^{-1}$ with that for~$k=0$ bounded even by~$(1+\frac{A'}{2d})^{-1}$. From \eqref{submartingale}, we thus get
\begin{equation}
%\label{}
\bigl|\psi(x)\bigr|^2\le\Bigl(1+\frac A{2d}\Bigr)^{2-2R}\Bigl(1+\frac{A'}{2d}\Bigr)^{-2}\,E^x\bigl(\1\{\tau\geq R\}|\psi(Y_{R})|^2\bigr),\qquad x\in D'.
\end{equation}
The reversibility of the simple random walk implies
\begin{equation}
%\label{}
\sum_{x\in D'} E^x\bigl(\1\{\tau\geq R\}|\psi(Y_{R})|^2\bigr)\le\sum_{x\in D'}\sum_{y\in D}P^x(Y_R=y)|\psi(y)|^2\le\Vert \psi\Vert_2^2,
\end{equation}
whereby the claim follows.
\end{proofsect}

\subsection{Rank-one perturbations}
To apply the above \emph{a priori} bounds in the proof of Theorem~\ref{thm-truncation}, we also need to tie this with rank-one perturbation arguments. First we prove a continuity statement:

\begin{lemma}
\label{lemma-D-U}
For $U\subset D$, let $\xi_s:=\xi-s\1_{D\smallsetminus U}$. For any $k\in\{1,\dots,|U|\}$, the map $s\mapsto\lambda_D^{\ssup k}(\xi_s)$ is non-increasing, Lipschitz continuous (with Lipschitz constant one) and
\begin{equation}
\label{D-U-limit}
\lambda_D^{\ssup k}(\xi_s)\,\underset{s\to\infty}\longrightarrow\,\lambda_U^{\ssup k}(\xi).
\end{equation}
\end{lemma}

\begin{proofsect}{Proof}
Since the eigenvalues are labeled in a decreasing order, the Minimax Theorem reads
\begin{equation}
\label{minimax}
\lambda_D^{\ssup k}(\xi)=\inf_{\begin{subarray}{c}
\HH_k\colon \HH_k\subset\C^{|D|}\\ \dim(\HH_k)=k-1
\end{subarray}}\,\,
\sup_{\begin{subarray}{c}
\phi\in \HH_k^{\perp}\\\|\phi\|_2=1
\end{subarray}}\,\,
\bigl(\phi,H_{D,\xi}\phi\bigr).
\end{equation}
The supremum goes over all $(k-1)$-dimensional linear subspaces of $\C^{|D|}$ and, here and henceforth, $(\psi,\phi):=\sum_{x\in D}\psi(x)^\star\phi(x)$ denotes the inner product in $\C^{|D|}$. From
\begin{equation}
\label{E:3.14}
H_{D,\xi_{s'}}=H_{D,\xi_s}+(s-s')\1_{D\smallsetminus U},
\end{equation}
we thus immediately get
\begin{equation}
%\label{}
0\le\lambda_D^{\ssup k}(\xi_{s'})-\lambda_D^{\ssup k}(\xi_s)\le s-s',\qquad 0\le s'<s.
\end{equation}
The Minimax Theorem also implies $\lambda_D^{\ssup k}(\xi_s)\ge\lambda_U^{\ssup k}(\xi)$ for all $k\in\{1,\dots,|U|\}$. Indeed, fix $\HH_k$ arbitrary with $\dim(\HH_k)=k-1$ but let $\phi\in\HH_k^\perp$ vanish outside $D\setminus U$. (Such a $\phi$ always exists because $|U|>k-1$.) Then $(\phi,H_{D,\xi}\phi)=(\phi,H_{U,\xi}\phi)\ge\lambda_U^{\ssup {\ell}}(\xi)$, where $\ell-1$ is the dimension of the projection of $\HH_k$ onto~$U$. Obviously $\ell\le k$ and so $(\phi,H_{D,\xi}\phi)\ge\lambda_U^{\ssup {k}}(\xi)$ as desired.

To prove the limit \eqref{D-U-limit}, consider a sequence $s_n\to\infty$ such that $\lambda_D^{\ssup k}(\xi_{s_n})\to\lambda_k$ as well as $\psi_{D,\xi_{s_n}}^{\ssup k}(x)\to\psi_k(x)$ exist. Since $\lambda_k$, with $k=1,\dots,|U|$, are finite by the above reasoning, we can immediately conclude that $(\Delta+\xi)\psi_k=\lambda_k\psi_k$ on~$U$ while $\psi_k(x)=0$ for $x\in D\setminus U$. As $k\mapsto\lambda_k$ is non-increasing, we must have $\lambda_U^{\ssup k}(\xi)=\lambda_k$.
\end{proofsect}

The use of Lemma~\ref{lemma-2-norm} will be aided by the following observation:

\begin{lemma}
\label{lemma-perturb}
Fix $A>0$, $R\in\N$ and, for $U$ as in Theorem~\ref{thm-truncation} and $s\ge0$, let $\xi_s:=\xi-s\1_{D\smallsetminus U}$.  Fix $k\in\{1,\dots,|D|\}$ and let $s>0$ be such that $\lambda_D^{\ssup k}(\xi_s)\ge\lambda_D^{\ssup 1}(\xi)-A$. Then
\begin{equation}
\label{E:3.12}
0\le \lambda_D^{\ssup k}(\xi)-\lambda_D^{\ssup k}(\xi_s)\le
\int_0^s\textd s'\Bigl(1+\frac A{2d}\Bigr)^{2-2R}\Bigl(1+\frac{A+s'}{2d}\Bigr)^{-2}.
\end{equation}
\end{lemma}

\begin{proofsect}{Proof}
The left inequality was claimed in Lemma~\ref{lemma-D-U}. The inequality on the right would be a consequence of the fact that $\frac\textd{\textd s}\lambda_D^{(k)}(\xi_s)=-\Vert\psi_{D,\xi_s}^{\ssup k}\1_{D\smallsetminus U}\Vert_2^2$ whenever the eigenvalue is non-degenerate. However, to get around the issue of degeneracy, we need to work a bit harder. 

Abbreviate 
\begin{equation}
%\label{}
g(s):=\Bigl(1+\frac A{2d}\Bigr)^{2-2R}\Bigl(1+\frac{A+s}{2d}\Bigr)^{-2}.
\end{equation}
We claim that it suffices to prove
\begin{equation}
\label{E:3.18b}
\lambda_D^{\ssup k}(\xi_s)\ge\lambda_D^{\ssup 1}(\xi)-A
\quad\Rightarrow\quad\limsup_{s'\uparrow s}\,\frac{\lambda_D^{\ssup k}(\xi_{s'})-\lambda_D^{\ssup k}(\xi_s)}{s-s'}
\le g(s).
\end{equation}
Indeed, if $f\colon\R\to\R$ is non-increasing and Lipschitz continuous, then $a<b$ implies
\begin{equation}
%\label{}
f(a)-f(b)=\lim_{h\downarrow0}\int_a^b\frac{f(x-h)-f(x)}h\,\textd x
\le\int_a^b\limsup_{h\downarrow0}\,\frac{f(x-h)-f(x)}h\,\textd x,
\end{equation}
by the Fatou Lemma. So if \eqref{E:3.18b} holds on an interval of $s$, the difference of the eigenvalues at the endpoints of this interval is bounded by the corresponding integral of~$g$.

To establish \eqref{E:3.18b}, consider a sequence $s_n'\uparrow s$ saturating the \emph{limes superior} and such that
\begin{equation}
%\label{}
\lambda_k:=\lim_{n\to\infty}\lambda_D^{\ssup k}(\xi_{s_n'})\quad\text{and}\quad
\psi_k(x):=\lim_{n\to\infty}\psi_{D,\xi_{s_n'}}^{\ssup k}(x)
\end{equation}
exist for all~$k\in\{1,\dots,|D|\}$ and all~$x\in D$. Obviously, $(\lambda_k,\psi_k)$ is an eigenvalue/eigenfunction pair for the field $\xi_s$ and $\lambda_k=\lambda_D^{\ssup k}(\xi_s)$. Both $\psi_k$ and $\psi_{D,\xi_{s'}}^{\ssup k}$ are eigenfunctions and so, by \eqref{E:3.14},
\begin{equation}
\begin{aligned}
\lambda_k\bigl(\psi_k,\psi_{D,\xi_{s'}}^{\ssup k}\bigr)
&=
\bigl(\psi_k,H_{D,\xi_s}\psi_{D,\xi_{s'}}^{\ssup k}\bigr)
\\
&=
\lambda_D^{\ssup k}(\xi_{s'})
\bigl(\psi_k,\psi_{D,\xi_{s'}}^{\ssup k}\bigr)+(s'-s)\bigl(\psi_k,\1_{D\smallsetminus U}\psi_{D,\xi_{s'}}^{\ssup k}\bigr).
\end{aligned}
\end{equation}
It follows that
\begin{equation}
%\label{}
\limsup_{s'\uparrow s}\frac{\lambda_D^{\ssup k}(\xi_{s'})-\lambda_D^{\ssup k}(\xi_s)}{s-s'}
=\frac{(\psi_k,\1_{D\smallsetminus U}\psi_k)}{(\psi_k,\psi_k)}
\end{equation}
and so, to get \eqref{E:3.18b}, it suffices to verify \eqref{E:3.8} for $\psi:=\psi_k$, $D':=D\setminus U$ and $A':=A+s$.

First note that $\lambda_D^{\ssup k}(\xi_s)\ge\lambda_D^{\ssup 1}(\xi)-A$ and $U\supset D_{R,A}(\xi)$ force $\xi<\lambda_D^{\ssup 1}(\xi)-2A\le\lambda_D^{\ssup k}(\xi_s)-A$ on $D\setminus U$. Hence, $\xi_s\le\lambda_k-A'$ on $D'$. Similarly, $\xi_s(x)\ge\lambda_k-A$ necessitates $\xi(x)\ge\lambda_D^{\ssup 1}(\xi)-2A$ which by $U\supset D_{R,A}(\xi)$ implies $\dist(x,D')\ge R$. The conditions of Lemma~\ref{lemma-2-norm} are thus met for the field~$\xi:=\xi_s$ and eigenvalue $\lambda:=\lambda_k$, and so \eqref{E:3.8} holds for $\psi:=\psi_k$ as desired.
\end{proofsect}

%\bf Something's wrong here. There is no restriction on $k$ and yet the result cannot hold for all $k=1,\dots,|D|$!!! \rm

The above facts are now assembled into the control of truncation from~$D$ to~$U$:

\begin{proofsect}{Proof of Theorem~\ref{thm-truncation}} 
Pick $A>0$ and $R\in\N$  such that the right-hand side of \eqref{eigen-distance} is $\le A/2$. Fix a $k\in\{1,\dots,|U|\}$ such that \eqref{lambda-theta-xi-initial} holds, i.e., $\lambda_D^{\ssup k}(\xi)\ge\lambda_D^{\ssup 1}(\xi)-\frac A2$. For~$s\ge0$, introduce the shifted field $\xi_s:=\xi-s\1_{D\smallsetminus U}$ and define
\begin{equation}
\label{tilde-s-def}
\tilde s:=\sup\bigl\{s\ge0\colon\lambda_D^{\ssup k}(\xi_s)\ge\lambda_D^{\ssup 1}(\xi)-A\bigr\}.
\end{equation}
By continuity of $s\mapsto\lambda_D^{\ssup k}(\xi_s)$, we have $\tilde s\in(0,\infty]$. Our aim is to show that $\tilde s=\infty$.

Suppose on the contrary that $\tilde s<\infty$. The bound \eqref{E:3.12} then shows that $\lambda_D^{\ssup k}(\xi)-\lambda_D^{\ssup k}(\xi_{\tilde s})$ is strictly less than the right-hand side of \eqref{eigen-distance} which by our assumption is~$\le\ffrac A2$. Therefore, 
\begin{equation}
%\label{}
\tilde s<\infty\quad\Longrightarrow\quad\lambda_D^{\ssup k}(\xi_{\tilde s})>\lambda_D^{\ssup k}(\xi)-\frac A2
\ge\lambda_D^{\ssup 1}(\xi)-A.
\end{equation}
This contradicts \eqref{tilde-s-def} because by continuity and monotonicity of $s\mapsto\lambda_D^{\ssup k}(\xi_s)$ we would have $\lambda_D^{\ssup k}(\xi_{s'})
\ge\lambda_D^{\ssup 1}(\xi)-A$ for an interval of $s'>\tilde s$. Hence $\tilde s=\infty$ as claimed.

To complete the proof we note that, in light of~\eqref{D-U-limit}, the difference $|\lambda_D^{\ssup k}(\xi)-\lambda_U^{\ssup k}(\xi)|$ is bounded by the integral in \eqref{E:3.12} with $s:=\infty$. A calculation then yields \eqref{eigen-distance}.
\end{proofsect}

\subsection{Spectral gap and potential profiles}
\label{sec-SpecGap}\noindent 
As an aside of the general arguments above, we will give a proof of Proposition~\ref{prop-eigengap} and establish versions of the inclusions \eqref{E:recall-union2} and \eqref{3:23-eq} that link local eigenvalues with potential profiles.

Let us first address Proposition~\ref{prop-eigengap}. A main tool of the proof is an inequality between the second eigenvalue and the principal eigenvalues in the nodal domains of the second eigenfunction. Note that $\psi_{C,\xi}^{\ssup 2}$, the second eigenfunction of $H_{C,\xi}$, has at least one negative and one positive value, since it is assumed orthogonal to the principal eigenfunction which is of one sign.

\begin{lemma}
%\label{lemma}
Let~$C\subset\Z^d$ be finite and define
\begin{equation}
\label{C-def}
B:=\bigl\{x\in C\colon \psi_{C,\xi}^{\ssup2}(x)\ge0\bigr\}.
\end{equation}
Then
\begin{equation}
%\label{}
\lambda_C^{\ssup 2}(\xi)\le\min\bigl\{\lambda_B^{\ssup 1}(\xi),\lambda_{C\smallsetminus B}^{\ssup 1}(\xi)\bigr\}.
\end{equation}
\end{lemma}

\begin{proofsect}{Proof}
Abbreviate $\psi:=\psi_{C,\xi}^{\ssup2}$. As~$\psi$ is orthogonal to the first eigenfunction, both~$B$ and~$C\setminus B$ are non-empty. The eigenvalue equation $(\Delta+\xi)\psi=\lambda_C^{\ssup 2}(\xi)\psi$ implies
\begin{equation}
\begin{aligned}
\label{vC-sum}
\lambda_C^{\ssup2}(\xi)\bigl\Vert \psi\1_B\bigr\Vert_2^2
&=\sum_{x\in B}\psi(x)(\Delta+\xi)\psi(x)
\\&= 
\sum_{x\in B}\Bigl(\,\xi(x)\psi(x)^2\,+\!\!\!\!\sum_{y\colon |y-x|=1}\psi(x)\bigl(\psi(y)-\psi(x)\bigr)\Bigr).\end{aligned}
\end{equation}
Let~$\hat \psi$ be equal to~$\psi$ on~$B$ and zero on~$C\setminus B$. Then
\begin{equation}
%\label{}
\psi(x)\bigl(\psi(y)-\psi(x)\bigr)\le \hat \psi(x)\bigl(\hat \psi(y)-\hat \psi(x)\bigr)
\end{equation}
for all pairs $x\in B$ and $y$ with~$|y-x|=1$. By the Minimax Theorem, the sum in \eqref{vC-sum} with~$\hat \psi$ instead of~$\psi$ is bounded by $\lambda_B^{\ssup 1}(\xi)\Vert\hat \psi\Vert_2^2$ and since $\Vert\hat \psi\Vert_2^2=\Vert \psi\1_B\Vert_2^2>0$, we thus get
\begin{equation}
%\label{}
\lambda_C^{\ssup 2}(\xi)\le\lambda_B^{\ssup 1}(\xi).
\end{equation}
The bound in terms of $\lambda_{C\smallsetminus B}^{\ssup 1}(\xi)$ is completely analogous.
\end{proofsect}

\begin{proofsect}{Proof of Proposition~\ref{prop-eigengap}}
Let~$B$ be as in \eqref{C-def}. As~$\LL_C(\xi)=\LL_B(\xi)+\LL_{C\smallsetminus B}(\xi)$, we may assume without loss of generality that~$\LL_B(\xi)\le\frac12\LL_C(\xi)$. Then
\begin{equation}
%\label{}
\lambda_C^{\ssup2}(\xi)-\rho\log\LL_C(\xi)
\le\lambda_B^{\ssup 1}(\xi)-\rho\log\LL_B(\xi)-\rho\log2
\le-\chi_B-\rho\log 2.
\end{equation}
Invoking \eqref{eigen-gap} and $\chi_B\ge\chi_C$ (as follows by a direct inspection of \eqref{chi-C-alter}) we get \eqref{eigen-chi-gap}.
\end{proofsect}

Next let us move to the inclusion \eqref{E:recall-union2}. We will in fact derive a stronger version by relating local eigenvalues to the quantity
\begin{equation}
%\label{}
\LL_{C,A}(\varphi):=\sum_{x\in C}\texte^{\varphi(x)/\rho}\1_{\{\varphi(x)\ge -2A\}},
\end{equation}
which is better suited for our later needs. Clearly, for $A:=\infty$ this degenerates to $\LL_C(\varphi)$. 

\begin{lemma}
\label{lemma-3.6q}
For all $a\in\R$, all finite~$C\subset\Z^d$ and all $A\ge\chi_C$ satisfying $A(1+\frac A{4d})\ge4d$,
\begin{equation}
\label{E:recall-union}
\bigl\{\xi\colon\lambda_C^{\ssup 1}(\xi)\ge a\bigr\}\subseteq\bigl\{\xi\colon\LL_{C,A}(\xi-a-\chi_C)\ge\texte^{-\eta(A)/\rho}\bigr\},
\end{equation}
where $\eta(A):=2d(1+\frac A{4d})^{-1}$.
\end{lemma}

\begin{proofsect}{Proof}
Let us first address the case $A:=\infty$. Suppose $\lambda_C^{\ssup 1}(\xi)\ge a$, let $\varphi:=\xi-a-\chi_C$ and note that $\lambda_C^{\ssup 1}(\varphi)\ge-\chi_C$.  We claim that $\LL_C(\varphi)\ge1$. Indeed, if we had $\LL_C(\varphi)<1$ then we could find an~$\epsilon>0$ such that $\tilde\varphi:=\varphi+\epsilon$ would obey $\LL_C(\tilde\varphi)\le1$ and yet $\lambda_C^{\ssup 1}(\tilde\varphi)>-\chi_C$. This would contradict \eqref{chi-C-alter}. Hence \eqref{E:recall-union} holds for $A:=\infty$.

Now take $A\ge\chi_C$, let $\varphi:=\xi-a-\chi_C$ and set $C':=\{x\in C\colon\varphi(x)\ge -2A\}$. As $\lambda_{C}^{\ssup 1}(\varphi)\ge-\chi_C$ implies $\varphi\le\lambda_{C}^{\ssup 1}(\varphi)-A$ on $C\setminus C'$, Theorem~\ref{thm-truncation} with $D:=C$, $U:=C'$ and~$R:=1$ shows --- thanks to the condition $A(1+\frac A{4d})\ge4d$ --- that $\lambda_{C'}^{\ssup 1}(\varphi)\ge\lambda_C^{\ssup1}(\varphi)-\eta(A)$. Therefore, $\lambda^{\ssup1}_C(\xi)\ge a$ gives $\lambda_{C'}^{\ssup 1}(\xi)\ge a-\eta(A)$ and so, by the claim for $A:=\infty$, $\LL_{C'}(\xi-a-\chi_{C'}+\eta(A))\ge1$. In conjunction with $\chi_{C'}\ge\chi_C$, this yields
\begin{equation}
\begin{aligned}
\LL_{C,A}(\xi-a-\chi_C)&=\LL_{C'}(\xi-a-\chi_C)
\\
&\ge\LL_{C'}\bigl(\xi-a-\chi_C'+\eta(A)\bigr)\texte^{-\eta(A)/\rho}\ge\texte^{-\eta(A)/\rho},
\end{aligned}
\end{equation}
as desired.
\end{proofsect}

As another aside we also recall that \eqref{chi-C-alter} implies that, for any finite~$C\subset\Z^d$, there is $\varphi_C\colon C\to\R$ such that $\LL_C(\varphi_C)\le1$ and $\lambda_{C}^{\ssup 1}(\varphi_C)=-\chi_C$.
Therefore,
\begin{equation}
\label{E:recall-lower}
\bigl\{\xi\colon\lambda_{C}^{\ssup 1}(\xi)\ge a\bigr\}\supseteq\bigl\{\xi\colon\xi\ge a+\chi_C+\varphi_C\bigr\}.
\end{equation}
This provides a bound that is in a sense opposite to \eqref{E:recall-union}.

Moving to the inclusion \eqref{3:23-eq}, similarly to \eqref{E:recall-union} we will restate this using the truncated functional $\LL_{C,A}$ as follows:

\begin{lemma}
\label{cor-gap}
For~$C\subset\Z^d$ finite, all~$a, a'\in\R$ and all $A\ge0$ sufficiently large,
\begin{equation}
\label{3:23-eq*}
\begin{aligned}
\lambda_C^{\ssup 1}(\xi)\ge a'&\,\,\text{ \rm\small AND }\,\,\lambda_C^{\ssup 1}(\xi)-\lambda_C^{\ssup 2}(\xi)\le \frac12\rho\log2\quad\Longrightarrow\quad\LL_{C,A}(\xi-a-\chi_C)\ge u,
\end{aligned}
\end{equation}
where $u$ is defined by
\begin{equation}
\label{3:23-eq2*}
\log u:=\frac{a'-a-\eta(A)}\rho+\frac12\log2
\end{equation}
with $\eta(A):=2d(1+\frac A{4d})^{-1}$.
\end{lemma}

\begin{proofsect}{Proof}
Let~$\xi$ be such that the conditions on the left of \eqref{3:23-eq} apply and, similarly as in the proof of Lemma~\ref{lemma-3.6q}, let $\varphi:=\xi-a-\chi_C$ and set $C':=\{x\in C\colon\varphi(x)\ge -2A\}$. For~$A$ large Theorem~\ref{thm-truncation} can be used; which then implies $|\lambda_C^{\ssup 2}(\xi)-\lambda_{C'}^{\ssup 2}(\xi)|\le\eta(A)$.  With the help of \eqref{eigen-chi-gap}, this yields
\begin{equation}
\begin{aligned}
a'\le\lambda_C^{\ssup 1}(\xi)
&\le \lambda_C^{\ssup 2}(\xi)+\frac12\rho\log2
\\
&\le \lambda_{C'}^{\ssup 2}(\xi)+\eta(A)+\frac12\rho\log2
\\
&\le\rho\log\LL_{C'}(\xi)-\chi_{C'}+\eta(A)-\frac12\rho\log2.
\end{aligned}
\end{equation}
Using $\chi_{C'}\ge\chi_C$ and $\LL_{C'}(\xi)=\texte^{(a+\chi_C)/\rho}\LL_{C,A}(\varphi)$, we get the claim.
\end{proofsect}

By a variant of estimates used in Lemma~\ref{lemma-3.6q}, we will now control the spatial concentration of the fields that are near optimizers of \eqref{chi-C-def}. This will be useful in the derivation of spatial decay of the corresponding eigenfunctions.

\begin{lemma}
\label{lemma-3.8w}
Define $A_0$ by $A_0(1+\frac{A_0}{4d})=4d$ and suppose~$A,\delta>0$ and $A':=-\frac12\rho\log(2\sinh\delta)$ satisfy $A\ge A'\ge d+A_0$ and $\eta(A)/\rho\le\delta$. There is $c=c(A,\delta)\in(0,\infty)$ such that for any $C\subset\Z^d$ finite, any $a\in\R$, any $r\ge1$ and any $\xi\colon C\to\R$ with 
\begin{equation}
\label{E:3.35ww}
\LL_{C,A}(\xi-a-\chi_C)\le\texte^\delta\quad\text{and}\quad
\lambda_C^{\ssup1}(\xi)\ge a+2d\Bigl(1+\frac{A'-d}{2d}\Bigr)^{1-2r},
\end{equation}
there is~$x\in C$ with the property
\begin{equation}
%\label{}
z\in C\,\,\,\&\,\,\,|z-x|\ge cr\quad\Rightarrow\quad\xi(z)\le\lambda_C^{\ssup1}(\varphi)-A'+\chi_C.
\end{equation}
\end{lemma}

\begin{proofsect}{Proof}
Abbreviate $\varphi(z):=\xi(z)-a-\chi_C$ and set $S:=\{z\in C\colon \varphi(z)>-2A'\}$. As $A'\le A$, every point of~$S$ contributes to~$\LL_{C,A}(\varphi)$. Hence, $\texte^\delta\ge\LL_{C,A}(\varphi)\ge\texte^{-2A/\rho}|S|$ and so $|S|\le\texte^{\delta+2A/\rho}$. Our goal is to use this to show that~$S$ also has a bounded diameter.

Given~$r\ge1$ as above, suppose, for the sake of contradiction, that~$S$ has diameter larger than~$2|S|r$. Then $S$ can be split into two  parts, $S=S_1\cup S_2$, such that $\dist(S_1,S_2)> 2r$. Pick~$x\in S_1$ and set $C':=\{z\in C\colon\dist(z,S_2)\le r\}$. Since $0<\dist(z,S_2)\le 2r$ implies
\begin{equation}
%\label{}
\xi(z)\le a+\chi_C-2A'\le\lambda_C^{\ssup1}(\xi)+\chi_C-2A'\le \lambda_C^{\ssup1}(\xi)-2(A'-d),
\end{equation}
and $A'-d\ge A_0$, we may use Theorem~\ref{thm-truncation} to conclude
\begin{equation}
%\label{}
\lambda_{C'}^{\ssup1}(\xi)\ge\lambda_C^{\ssup1}(\xi)-2d\Bigl(1+\frac{A'-d}{2d}\Bigr)^{1-2r}.
\end{equation}
By \eqref{E:3.35ww} the right-hand side is at least~$a$ and so Lemma~\ref{lemma-3.6q} gives $\LL_{C',A}(\varphi)\ge\texte^{-\eta(A)/\rho}\ge\texte^{-\delta}$, where the second bound comes from $\delta\ge\eta(A)/\rho$. But \eqref{E:3.35ww} and $x\in C\setminus C'$ also yield
\begin{equation}
\label{E:3.37ww}
\LL_{C',A}(\varphi)\le\LL_{C,A}(\varphi)-\texte^{\varphi(x)/\rho}
\le\texte^\delta-\texte^{\varphi(x)/\rho},
\end{equation}
whereby we get $\varphi(x)\le\rho\log(2\sinh\delta)=-2A'$. This contradicts~$x\in S$ and so $\diam S\le2|S|r$ must hold after all. Setting $c:=2\texte^{\delta+2A/\rho}$, the claim follows.
\end{proofsect} 

\begin{remark}
We note that considerable effort has been devoted to the study of the minimizers in the variational problem \eqref{chi-def}; cf G\"artner and den Hollander~\cite{GH99}. In spite of that, uniqueness of the minimizer remains open for small~$\rho$. In the same study, decay estimates for the eigenfunctions associated with the minimizing potential are provided.
\end{remark}

\section{Eigenvector localization: deterministic estimates}
\label{sec-EigenvLoc}\noindent 
Our discussion of deterministic estimates proceeds by giving the proof of Theorem~\ref{prop-exp-decay}. We will rely heavily on Theorem~\ref{thm-truncation} so, given $A>0$ and $R\in\N$, let us write
\begin{equation}
%\label{}
\epsilon_R:=2d\Bigl(1+\frac A{2d}\Bigr)^{1-2R}
\end{equation}
for the error bound in \eqref{eigen-distance}.
There are two main inputs into our proof of Theorem~\ref{prop-exp-decay}. The first of these is an inequality between the distance to the nearest eigenvalue and the ratio of masses that the eigenfunction puts on the boundary of a set relative to what it puts inside.

\begin{proposition}[Boundary mass vs gap]
\label{prop-localization}
For the setting of Theorem~\ref{prop-exp-decay}, suppose that~$U\subset D$ is such that
\begin{equation}
\label{U-separ}
z\in D\,\,\text{ \rm\small AND }\,\,\dist(z,\partial U)\leq R+1\quad\Longrightarrow\quad
\xi(z)\leq \lambda_D^{\ssup 1}(\xi)-2A.
\end{equation}
Put $U'=D\setminus(U\cup\partial U)$ (hence $\partial U=\partial U'$). Then
\begin{equation}
\label{prop4.9one}
\max\Big\{\frac{\|\psi\1_{\partial U}\|_2}{\|\psi\1_{U}\|_2},\frac{\|\psi\1_{\partial U'}\|_2}{\|\psi\1_{U'}\|_2}\Big\}\geq \frac {\gap(\lambda)-2\varepsilon_R}{8d}\wedge 1.
\end{equation}
\end{proposition}

The proof will require a specific comparison of the eigenvalues of $\Delta+\xi$ in~$U$ with different boundary conditions on~$\partial U$.

\begin{lemma}[Removing boundary condition]
\label{lemma-comparison}
Let~$U\subset\Z^d$ be finite and let~$\tilde \psi\colon\partial U\to\R$ obey for simplicity $\|\tilde \psi\1_{\partial U}\|_2\le 1$. Recall that~$\lambda_U^{\ssup k}(\xi)$ is the~$k$-th largest eigenvalue for operator~$\Delta+\xi$ in~$U$ with zero boundary condition, and let~$\tilde\lambda^{\ssup k}_U(\xi)$ be the~$k$-th largest eigenvalue for the same operator with boundary condition~$\tilde \psi$. Then
\begin{equation}
\label{4.41-odhad}
\bigl|\lambda_U^{\ssup k}(\xi)-\tilde \lambda_U^{\ssup k}(\xi)\bigr|\le 4d\|\tilde \psi\1_{\partial U}\|_2,
\qquad k=1,\dots,|U|.
\end{equation}
\end{lemma}

\begin{proofsect}{Proof}
Let~$\tilde\Delta$ denote the Laplace operator on~$U$ with boundary condition~$\tilde \psi$ on~$\partial U$. Then, for any  function $\psi\colon U\to \R$,
\begin{equation}
%\label{}
\begin{aligned}
\bigl(\psi,(\tilde\Delta+\xi)\psi\bigr)-\bigl(\psi,(\Delta+\xi)\psi\bigr)
%&=\sum_{\begin{subarray}{c}
%x\in U\\y\in \partial U\\|x-y|=1
%\end{subarray}}
%\Bigl(v(x)^2-\bigl[v(x)-\tilde v(y)\bigr]^2\Bigr)
=\sum_{\begin{subarray}{c}
x\in U,\,y\in \partial U\\|x-y|=1
\end{subarray}}\big(\tilde \psi(y)^2-2\psi(x)\tilde \psi(y)\big).
\end{aligned}
\end{equation}
Assuming $\Vert\psi\Vert_2=1$, the  Cauchy-Schwarz inequality and $\|\tilde \psi\1_{\partial U}\|_2\leq 1$ tell us
\begin{multline}
%\label{}
\qquad
\label{4.9one}
\Bigl|\bigl(\psi,(\tilde\Delta+\xi)\psi\bigr)-\bigl(\psi,(\Delta+\xi)\psi\bigr)\Bigr|
\\
\le 2d\|\tilde \psi\1_{\partial U}\|_2^2+2d\Vert \psi\Vert_2\|\tilde \psi\1_{\partial U}\|_2
\leq 4d\|\tilde \psi\1_{\partial U}\|_2.
\qquad
\end{multline}
By the Minimax Theorem (see \eqref{minimax}), $4d\|\tilde \psi\1_{\partial U}\|_2$ bounds the difference between the $k$-th largest eigenvalue of operators $\tilde\Delta+\xi$ and $\Delta+\xi$. Hence, \eqref{4.41-odhad} follows.
\end{proofsect}

\begin{proofsect}{Proof of Proposition~\ref{prop-localization}} 
Let $k\in\N$ be such that $\lambda=\lambda_D^{\ssup k}(\xi)$. Then the following facts hold:
\settowidth{\leftmargini}{(11)}
\begin{enumerate}
\item[(1)] 
By \eqref{U-separ}, the set $\partial U$ is at least~$R$ steps of the simple random walk from any point where $\xi$ exceeds $\lambda_D^{\ssup k}(\xi)-2A$. Since \eqref{E:3.43} holds, Theorem~\ref{thm-truncation} implies
\begin{equation}
\label{prop4.9two}
|\lambda_D^{\ssup k}(\xi)-\lambda_{D\smallsetminus\partial U}^{\ssup k}(\xi)|\leq \varepsilon_R.
\end{equation}
\item[(2)]
Restricting the eigenvalue relation on $D$ to $U$, resp., $U'$, there are $\ell_1,\ell_2\in\N$ such that
\begin{equation}
\label{kll'}
\lambda_D^{\ssup k}(\xi)=\tilde\lambda_U^{\ssup{\ell_1}}(\xi)=\tilde\lambda_{U'}^{\ssup {\ell_2}}(\xi),
\end{equation}
where $\tilde\lambda_U^{\ssup{\ell_1}}(\xi)$, resp., $\tilde\lambda_{U'}^{\ssup {\ell_2}}(\xi)$ is the $\ell_1$-th, resp., $\ell_2$-th eigenvalue of $\Delta+\xi$ in $U$, resp.,~$U'$ with boundary condition $\tilde \psi:=\psi/\|\psi\1_U\|_2$, resp., $\tilde \psi':=\psi/\|\psi\1_{U'}\|_2$ on $\partial U=\partial U'$.
\item[(3)]
If the left-hand side of \eqref{prop4.9one} is $\geq 1$, then there is nothing to prove, so let us assume that it is~$<1$. Lemma~\ref{lemma-comparison} then tells us
\begin{equation}
\label{E:ratios}
|\tilde\lambda_U^{\ssup {\ell_1}}(\xi)-\lambda_U^{\ssup {\ell_1}}(\xi)|\leq 4d \frac{\|\psi\1_{\partial U}\|_2}{\|\psi\1_{U}\|_2}
\end{equation}
and
\begin{equation}
\label{E:ratios2}
|\tilde\lambda_U^{\ssup {\ell_2}}(\xi)-\lambda_U^{\ssup {\ell_2}}(\xi)|\leq 4d \frac{\|\psi\1_{\partial U}\|_2}{\|\psi\1_{U'}\|_2}.
\end{equation}
\item[(4)] 
The Dirichlet eigenvalues in $D\setminus\partial U$ consist of the union of Dirichlet eigenvalues in~$U$ and~$U'$. It follows that there are $k_1,k_2\in\N$ such that
\begin{equation}
%\label{}
\lambda_U^{\ssup {\ell_1}}(\xi)=\lambda^{\ssup{k_1}}_{D\smallsetminus U}(\xi)
\quad\text{and}\quad
\lambda_{U'}^{\ssup {\ell_2}}(\xi)=\lambda^{\ssup{k_2}}_{D\smallsetminus U}(\xi).
\end{equation}
\end{enumerate}
Our goal is to show that $k\in\{k_1,k_2\}$ and $k_1\ne k_2$. 

First, Lemma~\ref{lemma-2-norm} tells us $\Vert \psi\1_{\partial U}\Vert_2\le\varepsilon_R/2d$ and so $\Vert \psi\1_U\Vert_2\vee\Vert \psi\1_{U'}\Vert_2\ge\frac12(1-\varepsilon_R/2d)\ge\ffrac14$. A calculation now shows that at least one of the right-hand sides in \twoeqref{E:ratios}{E:ratios2} is $\le8\varepsilon_R$ --- say the first one. But then
\begin{equation}
%\label{}
\bigl|\lambda_D^{\ssup k}(\xi)-\lambda_{D\smallsetminus U}^{\ssup{k_1}}(\xi)\bigr|=\bigl|\tilde\lambda_U^{\ssup {\ell_1}}(\xi)-\lambda_U^{\ssup {\ell_1}}(\xi)
\bigr|\le8\varepsilon_R<\gap(\lambda)-2\varepsilon_R.
\end{equation}
Since by \eqref{prop4.9two} $s\mapsto\lambda_D^{\ssup k}(\xi_s)$ stays at least $\gap(\lambda)-2\varepsilon_R$ away from other eigenvalues as $s$ slides off to infinity, 
%\bf NEED TO ASSUME \eqref{E:3.43} WITH $\epsilon_R$ ADDED ON RHS???\rm, 
we must have $k_1=k$. But if also $k_2=k$, then $H_{D\smallsetminus\partial U}$ would have two (degenerate) eigenvalues equal to $\lambda_{D\smallsetminus\partial U}^{\ssup k}(\xi)$ and so, by \eqref{prop4.9two}, $H_{D,\xi}$ would have another eigenvalue within $2\varepsilon_R$ of $\lambda_D^{\ssup k}(\xi)$. This is again impossible because $\gap(\lambda)>2\varepsilon_R$ and so $k_1\ne k_2$.

The rest of the proof now boils down to the estimates:
$$
\begin{aligned}
\gap(\lambda)-2\varepsilon_R
&\leq \bigl|\lambda_{D\smallsetminus\partial U}^{\ssup k}(\xi)-\lambda^{\ssup{k_2}}_{D\smallsetminus\partial U}(\xi)\bigr|
=\bigl|\lambda_U^{\ssup {\ell_1}}(\xi)-\lambda_U^{\ssup {\ell_2}}(\xi)\bigr|\\
&\leq \bigl |\tilde\lambda_U^{\ssup {\ell_1}}(\xi)-\lambda_U^{\ssup {\ell_1}}(\xi) \bigr |+ \bigl |\tilde\lambda_U^{\ssup {\ell_2}}(\xi)-\lambda_U^{\ssup {\ell_2}}(\xi) \bigr |\\
&\leq 4d \frac{\|\psi\1_{\partial U}\|_2}{\|\psi\1_{U}\|_2}+4d \frac{\|\psi\1_{\partial U}\|_2}{\|\psi\1_{U'}\|_2}\\
&\leq 8d\times \text{l.h.s.~of \eqref{prop4.9one}},
\end{aligned}
$$
whereby \eqref{prop4.9one} is finally proved.
\end{proofsect}

The second main input into the proof of Theorem~\ref{prop-exp-decay} is a continuity argument which we again state in general terms as follows: 

\begin{proposition}[Continuity argument]
\label{lemma-ind-step}
Fix $\delta\in(0,1)$ and $h>0$ so that the conditions expressed in \twoeqref{E:4.40}{E:4.41} hold. Fix $\d(x,\CC)$ satisfying (D0--D2) and suppose that a normalized eigenfunction~$\psi$ corresponding to eigenvalue~$\lambda$ of $H_{D,\xi}$ is such that, for some~$\CC\in\mathfrak C_{R,A}$,
\begin{equation}
\label{E:assume}
\bigl|\psi(z)\bigr|\le\texte^{-\delta h\d(z,\CC)},\qquad z\in D\setminus\CC,
\end{equation}
and, in addition, $\Vert \psi\1_\CC\Vert_2\ge\ffrac14$. Then we have $\Vert \psi\1_\CC\Vert_2>\ffrac12$ and
\begin{equation}
\label{E:conclusion}
\bigl|\psi(z)\bigr|<\texte^{-\delta h\d(z,\CC)},\qquad z\in D\setminus\CC.
\end{equation}
\end{proposition}

\begin{proofsect}{Proof}
Let us first consider~$z\in D\setminus D_{R,A}(\xi)$ and set $U:=B_R(z)\cap D$. For a path $(Y_j)_{j\in \N_0}$ of the simple random walk started at~$z$, let~$\tau_R$ denote the hitting time of $\partial U$. Recall the martingale~$M^\tau=(M_{\tau\wedge n})_{n\in\N_0}$ from Lemma~\ref{lemma-localization}; we will consider it at time $n:=\tau_R$. Note that if~$Y_{\tau\wedge\tau_R}\in\partial D$, then~$\psi(Y_{\tau_R})=0$; otherwise, the path $Y:=(Y_0,Y_1,\dots,Y_{\tau\wedge\tau_R-1})$ stays inside the set where $\xi<\lambda$. Since the loop-erasure of~$Y$ is an $R$-step self-avoiding nearest-neighbor path from $z$ to $\partial B_R(z)$, the product of the terms in \eqref{Mn-def} is at most~$\texte^{-hR}$, by \eqref{E:4.40}. (All the terms that we drop by loop-erasing are $\leq 1$.) The Optional Stopping Theorem implies
\begin{equation}
%\label{}
\bigl|\psi(z)\bigr|\le E^z(|M_{\tau\wedge\tau_R}|)\le\texte^{-hR}\max_{y\in \partial B_R(z)}\bigl|\psi(z)\bigr|.
\end{equation}
The assumption \eqref{E:assume} and (D1) allow us to bound the maximum by $\texte^{-\delta h[\d(z,\CC)-R]}$. (Note that we can seamlessly extend the maximum to all of $\partial B_R(z)$ because $\psi$ vanishes outside~$D$.)
Rearranging  terms we thus get
\begin{equation}
\label{E:3.55}
\bigl|\psi(z)\bigr|\le\texte^{-(1-\delta)hR}\texte^{-\delta h\d(z,\CC)},\qquad z\in D\setminus D_{R,A}(\xi).
\end{equation}
As~$1-\delta>0$ and $h>0$, this is even stronger than \eqref{E:conclusion}.

Our next goal is to boost the lower bound on $\Vert \psi\1_\CC\Vert_2$ from (non-strict) $\ffrac14$ to (strict) $\ffrac12$. To this end let $\kappa$ abbreviate the right-hand side of \eqref{prop4.9one} and note that, by \eqref{prop4.9one} with $U:=\CC$ and $U':=D\setminus(\CC\cup\partial\CC)$,
\begin{equation}
\label{E:3.56}
\min\bigl\{\Vert \psi\1_\CC\Vert_2,\Vert \psi\1_{D\smallsetminus(\CC\cup\partial\CC)}\Vert_2\bigr\}
\le\frac1\kappa\Vert\psi\1_{\partial\CC}\Vert_2.
\end{equation}
Since \eqref{E:3.55} applies to all $z\in\partial\CC$, we can estimate $\Vert\psi\1_{\partial\CC}\Vert_2$ as
\begin{equation}
%\label{}
\frac 1\kappa\Vert\psi\1_{\partial\CC}\Vert_2\le
\frac 1\kappa\sqrt{|\partial\CC'|}\,\max_{y\in\partial\CC}\bigl|\psi(y)\bigr|
\le \frac 1\kappa\sqrt{|\partial\CC|}\,\texte^{-(1-\delta)hR}.
\end{equation}
Invoking \eqref{E:4.41}, this is less than $\frac14$. By our assumption $\Vert \psi\1_\CC\Vert_2\ge\frac14$, we thus conclude from \eqref{E:3.56} that $\Vert \psi\1_{D\smallsetminus(\CC\cup\partial\CC)}\Vert_2\le\frac14$. But the same bound and $\kappa\le1$ yield
\begin{equation}
\Vert \psi\1_{D\smallsetminus \CC}\Vert_2^2
=\Vert \psi\1_{D\smallsetminus(\CC\cup\partial\CC)}\Vert_2^2+
\Vert\psi\1_{\partial\CC}\Vert_2^2\le\frac18,
\end{equation}
which implies $\Vert \psi\1_{\CC}\Vert_2\ge\sqrt{\,\ffrac78}>\ffrac12$, as desired.

We are now ready to prove \eqref{E:conclusion} for $z\in D_{R,A}(\xi)\setminus\CC$. Let~$\CC'$ denote the component that~$z$ belongs to. We will apply Proposition~\ref{prop-localization} with $U:=\CC'$ and $U':=D\setminus(\CC'\cup\partial\CC')$.  Since we already know $\Vert \psi\1_{\CC}\Vert_2>\ffrac12$,  we have $\Vert \psi\1_{U}\Vert_2<\Vert \psi\1_{U'}\Vert_2$ and so the maximum in \eqref{prop4.9one} is achieved by the term corresponding to~$U$. Therefore
\begin{equation}
%\label{}
\bigl|\psi(z)\bigr|\le\Vert \psi\1_{\CC'}\Vert_2\le\frac1\kappa\Vert \psi\1_{\partial\CC'}\Vert_2
\le\frac1\kappa\sqrt{|\partial\CC'|}\,\max_{y\in\partial\CC'}\bigl|\psi(y)\bigr|.
\end{equation}
By \eqref{E:3.55} we have $|\psi(y)|\le\texte^{-(1-\delta)hR-\delta h\d(y,\CC)}$ for each $y\in\partial\CC'$ and invoking condition~(D2) we can estimate $\d(y,\CC)\ge\d(z,\CC)-1$. Putting the terms together, we get
\begin{equation}
%\label{}
\bigl|\psi(z)\bigr|\le\Bigl(\,\frac1\kappa\sqrt{|\partial\CC'|}\texte^{\delta h-(1-\delta)hR}\Bigr)
\texte^{-\delta h\d(z,\CC)}.
\end{equation}
By \eqref{E:4.41} the prefactor in the large parentheses is~$\le\frac14$ for all~$\CC'$. The claim thus follows.
\end{proofsect}

\begin{proofsect}{Proof of Theorem~\ref{prop-exp-decay}}
We begin by identifying the component~$\CC$. Let~$k$ be such that $\lambda=\lambda_D^{\ssup k}(\xi)$, let~$U:=D\setminus D_{R,A}(\xi)$ and consider the deformation~$\xi_s:=\xi-s\1_{U}$ with $s\in(-\infty,\infty]$. By the assumption that $\gap(\lambda)>2\epsilon_R$, the eigenvalue $\lambda_D^{\ssup k}(\xi_s)$ stays non-degenerate for all $s\in[0,\infty]$, since this eigenvalue and its two neighbors $\lambda_D^{\ssup {k+1}}(\xi_s)$ and $\lambda_D^{\ssup {k-1}}(\xi_s)$ change by less than $\epsilon_R$ as~$s$ slides from $s:=0$ to $s:=\infty$. In other words, there is no eigenvalue crossing along the path.

It follows that also the corresponding eigenfunction $\psi_s$ changes continuously with~$s$ and its limit as~$s\to\infty$ exists and defines an eigenfunction $\bar \psi$ for~$\Delta+\xi$ in~$D\setminus U$ corresponding to~$\lambda_{D\smallsetminus U}^{\ssup k}(\xi)$. Clearly, there is a unique component of~$D\setminus U$ where~$\bar\psi$ puts all of its mass, because otherwise $\lambda_{D\smallsetminus U}^{\ssup k}(\xi)$ would be at least two-fold degenerate. We let~$\CC$ denote this component.

The bounds \eqref{E:assume} for $\psi$ replaced by $\psi_s$ and $\Vert \psi_s\1_\CC\Vert_2\ge\ffrac14$ are satisfied at $s:=\infty$; thanks to continuity in~$s$ they also hold for all~$s$ sufficiently large. Let~$s_0$ be the supremum of~$s\in\R$ where any of these bounds fails. We claim that~$s_0<0$. Indeed, (assuming $s_0>-\infty$) the bounds still hold (by continuity) at~$s=s_0$, and if~$s_0\ge0$, then Proposition~\ref{lemma-ind-step} would imply even the stronger bound $\Vert \psi_{s_0}\1_\CC\Vert_2\ge\ffrac12$. By continuity again, one could find~$\epsilon>0$ such that at least the bound $\Vert \psi_s\1_\CC\Vert_2\ge\ffrac14$ would hold for all~$s\in[s_0-\epsilon,s_0]$, in contradiction with the definition of~$s_0$. In particular, \eqref{E:assume} and $\Vert \psi\1_\CC\Vert_2\ge\ffrac14$ hold. As the eigenfunction is normalized and $\d(z,\CC)=0$ for $z\in\CC$, we have \eqref{eigenvdecay}, as desired.\end{proofsect}

\section{Coupling to i.i.d.\ random variables}
\label{sec-orderstat}
\noindent
Having dealt with the deterministic estimates that underly the proof of our main theorems, we now move on to the corresponding set of probabilistic arguments. Our specific task here is to establish a coupling to i.i.d.\ random variables as stated in Theorem~\ref{prop-4.6}. We assume throughout that $\xi=(\xi(x)\colon  x\in\Z^d)$ are i.i.d.\ random variables satisfying Assumption~\ref{ass}. 

\subsection{Extreme values of the fields}
\label{sec:extreme-val}\noindent
We begin by discussing the extreme order statistics of the potential field. Since $\xi(0)$ is continuously distributed, there is a unique $\hata_L$ such that
\begin{equation}
\label{E:4.8}
\BbbP(\xi(0)\ge\hata_L)=\frac1{L^d}.
\end{equation}
Assumption~\ref{ass} forces $F(s):=\log\log[\BbbP(\xi(0)>s)^{-1}]$ to grow, in the leading order, as a linear function with slope~$1/\rho$. Simple estimates then show
\begin{equation}
%\label{}
\hata_L=\bigl(\rho+o(1)\bigr)\log\log L,\qquad L\to\infty.
\end{equation}
The quantity $\hata_L$ plays the role of a centering sequence for the extreme order statistics of the field~$\xi$ in boxes of volume~$L^d$. Indeed, since $\texte^{F(\hata_L)}=d\log L$, setting 
\begin{equation}
%\label{}
b_L:=\frac\rho{d\log L}
\end{equation}
and using that ${F(\hata_L+sb_L)}={F(\hata_L)+\tfrac1\rho sb_L(1+o(1))}$ (by Assumption~\ref{ass}), we get
\begin{equation}
\label{maxxi}
\begin{aligned}
\BbbP\bigl(\,\max_{x\in B_L}\xi(x)\le\hata_L+sb_L\bigr)
&=\bigl(1-\exp\{-\texte^{F(\hata_L+sb_L)}\}\bigr)^{|B_L|}
\\
&=\exp\bigl\{-|B_L|\exp\{-\texte^{F(\hata_L+sb_L)}\}(1+o(1))\bigr\}
\\
&=\exp\Bigl\{-|B_L|\exp\bigl\{-\texte^{F(\hata_L)+\tfrac1\rho sb_L(1+o(1))}\bigr\}(1+o(1))\Bigr\}
\\
&\,\underset{L\to\infty}\longrightarrow\,\exp\bigl\{-\texte^{s}\bigr\}.
\end{aligned}
\end{equation}
It follows that the extreme-order law of $\xi$ lies in the Gumbel universality class. 

Our ultimate goal is to arrive at a similar conclusion also for the statistics of the maximal eigenvalues in~$D_L$. However, this requires more regularity than the above estimates, for two reasons: First, to make a (local) eigenvalue relevant requires arranging a whole ``profile'' of~$\xi$ values. Second, none of these $\xi$ values will reach into the extreme-order range --- i.e., within $O(b_L)$ of $\hata_L$ --- rather, they will all lie a positive deterministic constant below~$\hata_L$; cf~\eqref{E:1.16}. It turns out that we will need:

\begin{lemma}
\label{lemma-4.1}
Suppose Assumption~\ref{ass}. Then
\settowidth{\leftmargini}{(11)}
\begin{enumerate}
\item[(1)]
For any $\epsilon>0$ there is $r_0<\infty$ such that for all $r,r'$ with $r'\ge r\ge r_0$,
\begin{equation}
\label{E:asymp-1}
\texte^{(\frac1\rho-\epsilon)(r'-r)}\le
\frac{\log\BbbP(\xi(0)>r')}{\log\BbbP(\xi(0)>r)}\le\texte^{(\frac1\rho+\epsilon)(r'-r)}.
\end{equation}
\item[(2)]
The law of $\xi(0)$ has a density $f$ with respect to the Lebesgue measure which then satisfies
\begin{equation}\label{densityass}
\lim_{r\to\infty}\frac{f\big(r+s\texte^{-F(r)}\big)}{f(r)}=\texte^{-s/\rho}
\end{equation}
locally uniformly in~$s\in\R$.
\end{enumerate}
\end{lemma}

\begin{proofsect}{Proof}
For \eqref{E:asymp-1}, the ratio of the logs equals the exponential of $F(r')-F(r)$. Using the Fundamental Theorem of Calculus and Assumption~\ref{ass}, once $r',r$ are sufficiently large and $r'>r$,
\begin{equation}
\label{E:4.7}
\Bigl(\frac1\rho-\epsilon\Bigr)(r'-r)\le F(r')-F(r)\le\Bigl(\frac1\rho+\epsilon\Bigr)(r'-r).
\end{equation}
This implies part~(1) of the claim. For part (2), a calculation shows that
\begin{equation}
%\label{}
f(r):=F'(r)\,\BbbP\bigl(\xi(0)>r\bigr)\log\BbbP\bigl(\xi(0)>r\bigr)
\end{equation}
is the probability density of~$\xi(0)$, and so
\begin{equation}
%\label{}
\frac{f\big(r+s\texte^{-F(r)}\big)}{f(r)}
=\frac{F'(r+s\texte^{-F(r)})}{F'(r)}\,\frac{\BbbP\bigl(\xi(0)>r+s\texte^{-F(r)}\bigr)}{\BbbP\bigl(\xi(0)>r\bigr)}\,
\frac{\log\BbbP\bigl(\xi(0)>r+s\texte^{-F(r)}\bigr)}{\log\BbbP\bigl(\xi(0)>r\bigr)}.
\end{equation}
As~$\lim_{r\to\infty}F(r)=\infty$ for $r\to\infty$ by \eqref{E:4.2}, the first ratio on the right tends to one by \eqref{E:4.2} and the same applies to the third ratio by \eqref{E:4.7}. As for the middle ratio, we note
\begin{equation}
%\label{}
\frac{\BbbP\bigl(\xi(0)>r+s\texte^{-F(r)}\bigr)}{\BbbP\bigl(\xi(0)>r\bigr)}
=\exp\Bigl\{-\texte^{F(r)}\bigl(\texte^{F(r+s\texte^{-F(r)})-F(r)}-1\bigr)\Bigr\}.
\end{equation}
The claim follows by invoking  $F(r+s\texte^{-F(r)})-F(r)=\frac1\rho s\texte^{-F(r)}(1+o(1))$.
\end{proofsect}

With the help of \eqref{E:asymp-1} we now get the asymptotic formula \eqref{5.1-eq}.
This yields the identity in \eqref{E:4.11} which implies  that the only ``profiles'' of the field that we can realistically expect to see in~$D_L$ are those for which $\varphi:=\xi-\hata_L$ obeys $\LL(\varphi)\lesssim1$.

\subsection{Local eigenvalue estimates}
\label{sec-4.3}\noindent 
Many of our arguments that are to follow will require the following (rather crude) bounds for the principal eigenvalues in rectangular boxes of sublogarithmic size in~$L$. 

\begin{proposition}
\label{prop-4.2}
Let $R_L\to\infty$ with $R_L=o(\log L)$. Then for each~$\delta>0$ sufficiently small there is $\epsilon>0$ such that
\begin{equation}
\label{E:4.24q}
\BbbP\bigl(\lambda^{\ssup 1}_{B_{R_L}}(\xi)\ge \hata_L-\chi+\delta\bigr)\le L^{-d-\epsilon}
\end{equation}
and
\begin{equation}
\label{E:4.25q}
\BbbP\bigl(\lambda^{\ssup 1}_{B_{R_L}}(\xi)\ge \hata_L-\chi-\delta\bigr)\ge L^{-d+\epsilon}
\end{equation}
for all~$L$ sufficiently large. 
%In particular, if $a_L$ is defined as in \eqref{E:4.13q}, then
%\begin{equation}
%\label{E:4.26q}
%a_L=\hata_L-\chi+o(1),\qquad L\to\infty.
%\end{equation}
\end{proposition}

For the regime in \eqref{E:4.25q}, we will also need a statement with the quantifiers interchanged:

\begin{proposition}
\label{prop-4.2a}
Let $R_L\to\infty$ with $R_L=o(\log L)$. Then for each~$\epsilon>0$ there is $\delta>0$ such that
\begin{equation}
\label{E:4.26q}
\BbbP\bigl(\lambda^{\ssup 1}_{B_{R_L}}(\xi)\ge \hata_L-\chi-\delta\bigr)\le L^{-d+\epsilon}
\end{equation}
\end{proposition}

Let us start with the upper bounds \eqref{E:4.24q} and \eqref{E:4.26q}. Their proof will be based on the set inclusion \eqref{E:recall-union} established in Lemma~\ref{lemma-3.6q}. For that we will need to show that the event on the right of \eqref{E:recall-union} is dominated by configurations with nearly minimal value of~$\LL_{C,A}(\xi-a-\chi_C)$.

\begin{lemma}
\label{lemma-bigcup}
%Suppose Assumption~\ref{ass}. 
Let $\rho>0$ and $A>0$ be given. There is a constant $c_{A,\rho}$ such that for all $\alpha\in(0,\ffrac12)$, all $u\ge2/3$, all finite~$C\subset\Z^d$, all $L$ large enough and $d':=d2^{-2\alpha/3}\texte^{-A\alpha/(3\rho)}$,
\begin{equation}
\label{union}
\BbbP\bigl(\LL_{C,A}(\xi-\hata_L)\ge u\bigr)
\le c_{A,\rho}\biggl(\frac4\alpha\,\frac{A+\rho\log(4/3)}{\rho\log(4/3)}\biggr)^{|C|}\!\!L^{-d' u^{1-\alpha}}.
\end{equation}
\end{lemma}

\begin{proofsect}{Proof}
Fix $\rho>0$, $A>0$ and $\epsilon>0$. For $u>0$, define
\begin{equation}
%\label{}
F_u:=\bigl\{\xi\colon\LL_{C,A}(\xi-\hata_L)\in[u,2u)\bigr\}.
\end{equation}
The event in \eqref{union} is covered by $\bigcup_{n\ge0}F_{2^n u}$ so it suffices to derive a good estimate on~$\BbbP(F_u)$. Note that $2u>1$ because we assumed $u\ge2/3$.

Our first task is to derive a version of \eqref{5.1-eq} without the $o(1)$-term. For this, let~$r_0$ be such that \eqref{E:asymp-1} holds. We will also assume that $L$ is so large that $\hata_L-A\ge r_0$. We claim that then
\begin{equation}
\label{E:4.17}
-A\le s\le\rho\log(2u)\quad\Rightarrow\quad\BbbP\bigl(\xi(0)\ge \hata_L+s\bigr)\le
\bigl(L^{-d (2u)^{-2\rho\epsilon}\texte^{-A\epsilon}}\bigr)^{\texte^{s/\rho}}.
\end{equation}
Indeed, let~$\rho'$ be defined by $\ffrac1{\rho'}:=\ffrac1\rho+\epsilon$. Then \eqref{E:asymp-1} bounds the probability by $L^{-d\theta}$ with $\theta:=\texte^{s/\rho-\epsilon|s|}$. Since  $s\le\rho\log(2u)$ implies $|s|\le 2\rho\log(2u)-s$, this is further bounded by $L^{-d\theta'}$ where $\theta':=(2u)^{-2\rho\epsilon}\texte^{s/\rho'}$. As $\texte^{s/\rho'}\ge \texte^{-A\epsilon+s/\rho}$, \eqref{E:4.17} follows.

Next we will discretize the set of possible potential values to cover~$F_u$ by a finite union of sets. Set $\delta:=\epsilon\rho^2\log(2u)$ and let $S_{A,u}:=\{-A+m\delta\colon m\in\N_0\}\cap(-\infty,\rho\log(2u)]$. Define
\begin{equation}
%\label{}
\varphi_\xi(x):=-A+\delta\bigl\lfloor(\xi(x)+A)/\delta\bigr\rfloor.
\end{equation}
Then $\varphi$ takes values in $S_{A,u}$ and 
\begin{equation}
\label{E:5.22q}
\LL_{C,A}(\xi-\hata_L)\ge u\quad\Rightarrow\quad\LL_{C,A}(\varphi_\xi)\ge u\texte^{-\delta/\rho}\ge u^{1-\rho\epsilon}.
\end{equation}
Let~$\{C_k\colon k=1,\dots,2^{|C|}-1\}$ be an enumeration of all non-empty subsets of~$C$ and let $\{\varphi_{k,j}\}$ be an enumeration of all functions $\varphi_{k,j}\colon C_k\to S_{A,u}$. \eqref{E:5.22q} then gives
\begin{equation}
\label{}
F_u\,\subseteq\,\,\bigcup_k\,\bigcup_{j\colon\LL_{C_k}(\varphi_{k,j})\ge u\texte^{-\delta/\rho}}\bigl\{\xi\colon\xi\ge\hata_L+\varphi_{k,j}\text{ on }C_k\bigr\}.
\end{equation}
Denoting $d'':=d2^{-2\rho\epsilon}\texte^{-A\epsilon}$, the condition $\LL_{C,A}(\varphi_{k,j})\ge u^{1-\rho\epsilon}$ implies
\begin{equation}
%\label{}
\BbbP\bigl(\xi\colon\xi\ge\hata_L+\varphi_{k,j}\text{ on }C_k\bigr)\le L^{-d''u^{1-3\rho\epsilon}}.
\end{equation}
The total number of pairs $(k,j)$ contributing to the union is at most
\begin{equation}
%\label{}
\Bigl(1+\frac{A+\rho\log(2u)}\delta\Bigr)^{|C|}=\biggl(1+\,\frac{A+\rho\log(4/3)}{\rho^2\epsilon\log(4/3)}\biggr)^{|C|}.
\end{equation}
So if~$\epsilon$ (so far arbitrary) is linked to $\alpha$ via $\alpha:=3\rho\epsilon$, we get
\begin{equation}
\label{E:5.21y}
\BbbP(F_u)\le \biggl(\frac4\alpha\,\frac{A+\rho\log(4/3)}{\rho\log(4/3)}\biggr)^{|C|}\!\!L^{-d' u^{1-\alpha}}
\end{equation}
with $d'$ related to~$\alpha$ as in the claim.
Since the prefactor is independent of~$u$, the desired bound follows by summing \eqref{E:5.21y} over $u$ taking values in $\{2^nu\colon n\ge0\}$. (Note that the restriction on~$L$ was independent of~$u$ and $\alpha$.)
\end{proofsect}

We can now prove the upper bounds in Propositions~\ref{prop-4.2} and \ref{prop-4.2a}:

\begin{proofsect}{Proof of \eqref{E:4.24q} and \eqref{E:4.26q}} 
For $A\ge\chi_C$ and $\delta\in\R$, the inclusion \eqref{E:recall-union} and $\chi_C\ge\chi$ show
\begin{equation}
%\label{}
\bigl\{\lambda^{\ssup1}_C(\xi)\ge\hata_L-\chi+\delta\bigr\}\subseteq \bigl\{\xi\colon\LL_{C,A}(\xi-\hata_L)\ge u\bigr\}\quad\text{for}\quad u:=\texte^{(\delta-\eta(A))/\rho}.
\end{equation}
Applying \eqref{union} with $C:=B_{R_L}$, and noting that the term exponential in~$|C|$ is $L^{o(1)}$, the desired probability is at most $L^{d' u^{1-\alpha}+o(1)}$, for any $\alpha\in(0,\ffrac12)$, provided $\delta$ is small and $A$ large so that $u\ge2/3$. Now for~$A$ large and $\alpha>0$ small, $d-d'$ can be made arbitrary small (positive) while~$u$, which satisfies $u>1$ for $\delta>0$ and $u<1$ for $\delta<0$, can be made as close to one as desired by choosing $\delta$ small. This proves the desired bounds.
\end{proofsect}

Concerning the lower bound in \eqref{E:4.25q}, we first state:

\begin{lemma}
\label{lemma5.5}
For each~$\rho'>\rho$ there is $K_{\rho'}<\infty$ such that for any finite $C\subset\Z^d$, any $\delta\in[0,1)$ and any~$L$, we have
\begin{equation}
\label{E:4.27q}
\BbbP\bigl(\lambda_{C}^{\ssup 1}\ge\hata_L-\chi_C-\delta\bigr)\ge L^{-d\theta}
\quad\text{for}\quad
\theta:=\texte^{-\delta/\rho'}\bigl(1+K_\rho'|C|^{\rho/\rho'}\texte^{-\hata_L/\rho'}\bigr).
\end{equation}
\end{lemma}

\begin{proofsect}{Proof}
Recall \eqref{E:recall-lower} and the notation used therein. Fix~$\rho'>\rho$ and let $r_0$ be such that the bound on the right \eqref{E:asymp-1} holds for all $r,r'\ge r_0$ with $1/\rho-\epsilon$ equal to~$1/\rho'$. Setting $\varphi_C':=\varphi_C\vee(r_0-\hata_L+\delta)$ and noting that $\varphi_C'\ge\varphi_C$, we get
\begin{equation}
\begin{aligned}
\BbbP\bigl(\lambda_{C}^{\ssup 1}\ge\hata_L-\chi_C-\delta\bigr)
&\ge\BbbP\bigl(\xi\ge\hata_L-\delta+\varphi_C'\bigr)
\\
&\ge\exp\Bigl\{-\log(L^d)\texte^{-\delta/\rho'}\sum_{x\in C}\texte^{\varphi_C'(x)/\rho'}\Bigr\}
\end{aligned}
\end{equation}
But $\rho'>\rho$ implies $\sum_{x\in C}\texte^{\varphi(x)/\rho'}\le\LL_C(\varphi)^{\rho/\rho'}$ and for $\LL_C(\varphi_C')$ we get
\begin{equation}
%\label{}
\LL_C(\varphi_C')\le\LL_C(\varphi_C)+|C|\texte^{(r_0+\delta-\hata_L)/\rho}.
\end{equation}
In light of $\LL_C(\varphi_C)\le1$ and a simple convexity estimate, we get \eqref{E:4.27q} with $K_{\rho'}:=\texte^{(r_0+1)/\rho'}$.
\end{proofsect}

\begin{proofsect}{Proof of \eqref{E:4.25q}}
Since $\chi_{B_{R_L}}\downarrow\chi$ and $|B_{R_L}|=(\log L)^{o(1)}$ while $\hata_L\asymp\log\log L$, the bound \eqref{E:4.27q} with $C:=B_{R_L}$ yields \eqref{E:4.25q} for $L$ large as soon as $d-\epsilon> d\texte^{-\delta/\rho}$.
\end{proofsect}

%\bf MORE TO FOLLOW ... \rm

%\begin{comment}

\subsection{Approximation by i.i.d.\ process}
\label{sec-highpeaks}\noindent 
We are now ready to assemble the arguments needed for the proof of Theorem~\ref{prop-4.6}. Let $R_L$ be a sequence of integers subject to \eqref{E:RLcond}. As before, let~$\mathfrak C_{R_L,A}$ denote the set of connected components of $D_{R_L,A}(\xi)$ for $D:=D_L$.

\begin{lemma}
\label{lemma-4.7a}
For any~$A>0$ there is an integer $n_A<\infty$ such that
\begin{equation}
%\label{}
\diam\eusm C\le n_A R_L,\qquad\forall\eusm C\in\mathfrak C_{R_L,A},
\end{equation}
occurs with probability tending to one as~$L\to\infty$.
\end{lemma}

\begin{proofsect}{Proof}
This is a consequence of \eqref{5.1-eq} and a straightforward union bound. Indeed, let $F_{L,n}(x)$ denote the event that $B_{nR_L}(x)$ contains at least~$n$ vertices~$z$ with $\xi(z)\ge\hata_L-\chi-2A$. Then
\begin{equation}
\label{E:4.37q}
\bigl\{\exists\eusm C\in\mathfrak C_{R_L,A}\colon\diam\eusm C>nR_L\bigr\}\subseteq\bigcup_{x\in B_{nL}(0)}F_{L,n}(x).
\end{equation}
By \eqref{5.1-eq} and a union bound we obtain
\begin{equation}
%\label{}
\BbbP\bigl(F_{L,n}(x)\bigr)\le \bigl|B_{nr_L}\bigr|^n L^{-dn\theta+o(1)}\quad\text{where}\quad\theta:=\texte^{-(\chi+2A)/\rho}.
\end{equation}
Since $R_L=L^{o(1)}$, as soon as $n$ is so large that $n\theta>1$, the probability of the union on the right of \eqref{E:4.37q} will tend to zero as~$L\to\infty$.
\end{proofsect}

Next we will focus attention on components where the eigenvalue is close to the optimal threshold $\hata_L-\chi$. For these components we get:

\begin{lemma}
\label{lemma-4.7}
Given~$A>0$ sufficiently large, there is $\delta>0$ such that the following holds with probability tending to one as $L\to\infty$: For any $\eusm C\in\mathfrak C_{R_L,A}$ that obeys $\lambda^{\ssup 1}_{\eusm C}(\xi)\ge\hata_L-\chi-\delta$,
\settowidth{\leftmargini}{(1111)}
\begin{enumerate}
\item[(1)] $\lambda^{\ssup 1}_{\eusm C}(\xi)-\lambda^{\ssup 2}_{\eusm C}(\xi)\ge\frac12\rho\log2$,
\item[(2)] $\eusm C\subset B_{N_L}(x)$ for some $x\in\bigl((N_L+1)\Z)^d$ and, in addition, $\dist(\eusm C,D_L^\cc)>N_L$. 
\end{enumerate}
If, in addition, $\lambda^{\ssup 1}_{\eusm C'}(\xi)\ge\hata_L-\chi-\delta$ for some $\eusm C'\in\mathfrak C_{R_L,A}$ then 
%\settowidth{\leftmargini}{(11)}
\begin{enumerate}
\item[(3)]
either $\eusm C'=\eusm C$ or $\dist(\eusm C',\eusm C)>N_L$.
\end{enumerate}
\end{lemma}

\begin{proofsect}{Proof}
Assume that $A$ is large and pick $\eusm C\in\mathfrak C_{R_L,A}$ that obeys $\lambda^{\ssup 1}_{\eusm C}(\xi)\ge\hata_L-\chi-\delta$. By Lemma~\ref{lemma-4.7a}, we may also assume $\diam\eusm C\le n_A R_L$. For the claim in~(1), if we had $\lambda^{\ssup 1}_{\eusm C}(\xi)-\lambda^{\ssup 2}_{\eusm C}(\xi)\le\frac12\rho\log2$, then Lemma~\ref{cor-gap} with the choices
\begin{equation}
%\label{}
a:=\hata_L-\chi_{\eusm C}\quad\text{and}\quad a':=\hata_L-\chi-\delta
\end{equation}
would yield $\LL_{\eusm C,A}(\xi-\hata_L)\ge\sqrt2\texte^{-\eta(A)/\rho}=:u$. By Lemma~\ref{lemma-bigcup} and the fact that $|\eusm C|=O(R_L^d)=o(\log L)$, the probability that a given set~$C$ is a component with these properties is at most $L^{-d\sqrt u+o(1)}$. But there are at most $L^d\texte^{O(R^d)}=L^{d+o(1)}$ ways to choose such a connected component in~$D_L$ and so~(1) follows by a union bound and $u>1$.

For (2) and (3), let us abbreviate $r_L:=n_A R_L$. Given $\epsilon>0$ let $\delta>0$ be as in Proposition~\ref{prop-4.2}. Since $\eusm C\subset B_{r_L}$ implies $\lambda^{\ssup 1}_{B_{r_L}}(\xi)\ge\lambda^{\ssup 1}_{\eusm C}(\xi)$, from \eqref{E:4.26q} we immediately have
\begin{equation}
%\label{}
\BbbP\Bigl(\exists\eusm C\in\mathfrak C_{R_L,A}\colon x\in\eusm C,\,\diam(\eusm C)\le r_L,\,\lambda^{\ssup 1}_{\eusm C}(\xi)\ge \hata_L-\chi-\delta\Bigr)\leq |B_{r_L}|L^{-d+\epsilon}
\end{equation}
for any $x\in B_{nL}(0)$. Now if (2) fails for some component $\eusm C\in\mathfrak C_{R_L,A}$, then $\eusm C$ contains a vertex either in $D_L\setminus\bigcup_{x\in((N_L+1)\Z)^d}B_{N_L}(x)$ or in $\{x\in D_L\colon\dist(x,D_L^\cc)\le N_L\}$. The former set has cardinality of order $L^d N_L^{-1}$ while the latter has cardinality of order $N_L L^{d-1}$ (indeed, thanks to recifiability of~$\partial D$ we have $\partial D_L=O(L^{d-1})$). Hence, for some constant $c_1$,
\begin{equation}
%\label{}
\BbbP\bigl(\text{(2) fails}\bigr)\le c_1\bigl(L^d N_L^{-1}+N_L L^{d-1}\bigr)R_L^d L^{-d+\epsilon},
\end{equation}
For (3) a similar argument yields
\begin{equation}
%\label{}
\BbbP(\bigl(\text{(3) fails}\bigr)\le c_2 L^d N_L^d (R_L^{d}L^{-d+\epsilon})^2.
\end{equation}
Using \twoeqref{E:RLcond}{E:RLNL}, both of these tend to zero as~$L\to\infty$ once $\epsilon$ is small enough (but fixed).
\end{proofsect}

Finally, we also need a (slightly more explicit) version of \eqref{E:1.16}:

\begin{lemma}
\label{lemma-4.9}
We have
\begin{equation}
\label{E:5.31s}
\lambda_{D_L}^{\ssup1}(\xi)-\hata_L\,\underset{L\to\infty}\longrightarrow\,-\chi,
\qquad\text{\rm in $\BbbP$-probability}.
\end{equation}
\end{lemma}

\begin{proofsect}{Proof}
Fix $\epsilon>0$ and consider the event that $\lambda_{D_L}^{\ssup1}(\xi)\ge\hata_L-\chi-\epsilon$.
Cover~$D_L$ by order~$(L/R_L)^d$ disjoint translates of~$B_{R_L}$. By \eqref{E:4.25q} and $R_L=L^{o(1)}$, with probability tending to one as~$L\to\infty$, at least in one of these boxes the principal Dirichlet eigenvalue exceeds $\hata_L-\chi-\delta$. Since $\lambda_{D_L}^{\ssup1}$ dominates all these eigenvalues, we get a lower bound in \eqref{E:5.31s}.

Next let us examine the event $F:=\{\lambda_{D_L}^{\ssup1}(\xi)\ge\hata_L-\chi+\epsilon\}$. Let $A>0$ and fix $R_L\to\infty$ with $R_L=o(\log L)$. Assume, with the help of Lemma~\ref{lemma-4.7a}, that all components of $\mathfrak C_{R_L,A}$ have diameter less than $n_A R_L$. Theorem~\ref{thm-truncation} thus implies
\begin{equation}
%\label{}
\lambda_{D_L}^{\ssup1}(\xi)\le\max_{x\in D_L}\lambda_{B_{n_A R_L}(x)}^{\ssup1}(\xi)+\frac\epsilon2
\end{equation}
with probability tending to one as $L\to\infty$. So on~$F$, at least one of the boxes $B_{n_A R_L}(x)$, with $x\in D_L$, has $\lambda_{B_{n_A R_L}(x)}^{\ssup1}(\xi)\ge\hata_L+\ffrac\epsilon2$. By Proposition~\ref{prop-4.2}, this has probability $o(L^{-d})$ and, since $D_L=O(L^d)$, also an upper bound in \eqref{E:5.31s} holds.
\end{proofsect}

We are now finally ready to establish the coupling of the top part of the spectrum in~$D_L$ to a collection of i.i.d.\ random variables.

\begin{proofsect}{Proof of Theorem~\ref{prop-4.6}}
Fix $A>0$ large and let $\delta>0$ be such that the conclusions of Lemmas~\ref{lemma-4.7a},~\ref{lemma-4.7} and~\ref{lemma-4.9} hold. Let~$U$ denote the set $D_{R_L,A}$ for~$D:=D_L$. Take $\delta'<\min\{\delta/2,A/2\}$. Theorem~\ref{thm-truncation} then implies, for all $k=1,\dots,|D_L|$,
\begin{equation}
\label{E:prvni}
\lambda_{D_L}^{\ssup k}(\xi)\ge\hata_L-\chi-\delta'\quad\Rightarrow\quad
\bigl|\lambda_{D_L}^{\ssup k}(\xi)-\lambda_{U}^{\ssup k}(\xi)\bigr|\le 2d\Bigl(1+\frac A{2d}\Bigr)^{1-2R_L}.
\end{equation}
But the spectrum in~$U$ is the union of the spectra in the components in~$\mathfrak C_{R_L,A}$ and, once conclusions (2-3) in Lemma~\ref{lemma-4.7} are in place, we only need to pay attention to components that are entirely contained, and single of that kind, in one of the boxes $B_{N_L}^{(i)}$, $i=1,\dots, m_L$. Since Lemma~\ref{lemma-4.7}(1) tells us that we can also disregard all but the principal eigenvalue, if $\eusm C\subset B_{N_L}^{(i)}$ is such a component, Theorem~\ref{thm-truncation} yields
\begin{equation}
\label{E:druhy}
\bigl|\lambda_{\eusm C}^{\ssup 1}(\xi)-\lambda_k(\xi)\bigr|\le 2d\Bigl(1+\frac A{2d}\Bigr)^{1-2R_L}.
\end{equation}
Combining \twoeqref{E:prvni}{E:druhy}, the claim follows.
\end{proofsect}

\section{Eigenvalue order statistics}
\label{sec6}\noindent
Our next item of business is a proof of extreme order statistics for eigenvalues in $D_L$ as $L\to\infty$. Having coupled the eigenvalues at the top of the spectrum of $H_{D_L,\xi}$ to a collection of i.i.d.\ random variables --- namely the principal eigenvalues in disjoint subboxes of side~$N_L$ --- the argument is reduced to identifying the max-order class that these variables fall into.

\subsection{Determining the max-order class}
\label{sec-localeigenv}\noindent
Our strategy is to first identify the max-order class for eigenvalues in boxes of side~$R_L$ and only then relate this to the eigenvalues in boxes of side~$N_L$.

\begin{proposition}%[Max-order class of local eigenvalues]
\label{prop-smallbox}
Suppose Assumption~\ref{ass} holds and, for $R_L\to\infty$ with $R_L=(\log L)^{o(1)}$,
let $a_L$ be as defined in \eqref{E:4.13q}. Let~$b_L$ obey
\begin{equation}
\label{E:4.46q}
b_L\log L\,\underset{\,n\to\infty\,}\longrightarrow\,\frac \rho d.
\end{equation}
Then, with $o(1)\to0$ as~$L\to\infty$ uniformly on compact set of~$s,r\in\R$, 
\begin{equation}
\label{4.16-eq}
\BbbP\bigl(\lambda_{B_{R_L}}^{\ssup 1}\ge a_L+rb_L\bigr)=\texte^{-s+o(1)}\,\BbbP\bigl(\lambda_{B_{R_L}}^{\ssup 1}\ge a_L+(r-s)b_L\bigr).
\end{equation}
\end{proposition}

\begin{remark}
It is the proof of this proposition that requires us to assume that the law of~$\xi(0)$ has a density with respect to the Lebesgue measure. Although this restriction can be overcome to some extent, we have not succeeded in formulating a more general condition that would yield a comparably easy proof of the asymptotic \eqref{4.16-eq}. A natural idea how to deal with discontinuous laws would be to first approximate the spectrum by that of a continuously-distributed field and then apply the present approach.
\end{remark}

The main idea of the proof of Proposition~\ref{prop-smallbox} is to compensate for a shift in the eigenvalue by way of a rigid shift of the field configuration. In order to keep the action confined to the asymptotic regime, we will only shift the values of $\xi$ that are close to~$\hata_L$. Given $A>0$ and $L\ge1$, consider the continuous function $g_{L,A}\colon\R\times\R\to\R$ given by
\begin{equation}
g_{L,A}(\xi,s):=\begin{cases}
\xi-s,\qquad&\text{if }\xi\ge\hata_L-A,
\\
\xi,\qquad&\text{if }\xi\le\hata_L-2A,
\\
\text{linear},\qquad&\text{else}.
\end{cases}
\end{equation}
Clearly, for $s<A$, the map $\xi\mapsto g_{L,A}(\xi,s)$ is strictly increasing. The deterministic part of the change-of-measure argument is provided by:

\begin{lemma}
%\label{lemma}
Given a finite~$C\subset\Z^d$, a configuration~$(\xi(x))_{x\in C}$ and $A>0$, abbreviate
\begin{equation}
\label{E:4.48w}
\tilde\xi_s(x):=g_{L,A}\bigl(\xi(x),s\bigr).
\end{equation}
Then for all $a\in\R$ and $s\ge0$,
\begin{equation}
\label{E:4.49w}
\bigl\{\lambda_C^{\ssup1}(\xi)\ge a+s\bigr\}\subseteq\bigl\{\lambda_C^{\ssup1}(\tilde\xi_{s})\ge a\bigr\}
\end{equation}
and, for all $a\ge \hata_L-A/2$ and all $s\ge0$,
\begin{equation}
\label{E:4.50w}
\bigl\{\lambda_C^{\ssup1}(\xi)\ge a+s\bigr\}\supseteq\bigl\{\lambda_C^{\ssup1}(\tilde\xi_{s'})\ge a\bigr\},
\end{equation}
where $s':=s/[1-2d(1+\frac A{4d})^{-2}]$.
\end{lemma}

\begin{proofsect}{Proof}
Abbreviate $O:=\{x\in C\colon \xi(x)\ge \hata_L-A\}$ and note that $\tilde\xi_s=\xi(x)-s$ for~$x\in O$.
Since~$s\ge0$, the variational characterization of the principal eigenvalue tells us
\begin{equation}
\label{E:4.51w}
\lambda_C^{\ssup1}(\tilde\xi_{s})+s\ge \lambda_C^{\ssup1}(\xi)\ge\lambda_C^{\ssup1}(\tilde\xi_{s})+s\sum_{x\in O}\bigl|\psi_{C,\tilde\xi_{s}}^{\ssup1}(x)\bigr|^2.
\end{equation}
The inequality on the left then immediately yields \eqref{E:4.49w}. 

For \eqref{E:4.50w}, let $s'$ be as given and let us assume $\lambda_C^{\ssup1}(\tilde\xi_{s'})\ge a$. Then $a\ge\hata_L-A/2$ and $\tilde\xi_{s'}\le \hata_L-A$ on $C\setminus O$ imply $\tilde\xi_{s'}\le a-A/2\le\lambda_C^{\ssup1}(\tilde\xi_{s'})-A/2$ on $C\setminus O$ and thus by Lemma~\ref{lemma-2-norm} with $D':=C\setminus O$, $A':=A$ and $R:=1$ and $A$ replaced by $A/2$,
\begin{equation}
%\label{}
\sum_{x\in O}\bigl|\psi_{C,\tilde\xi_{s'}}^{\ssup1}(x)\bigr|^2\ge 1-2d\Bigl(1+\frac A{4d}\Bigr)^{-2}.
\end{equation}
The inequality on the right of \eqref{E:4.51w} with~$s$ replaced by~$s'$ yields $\lambda_C^{\ssup1}(\xi)\ge\lambda_C^{\ssup1}(\tilde\xi_{s'})+s\ge a+s$.
\end{proofsect}

The shift of the field will give rise to a term reflecting the change in the underlying measure. This term can be evaluated rather explicitly. As already pointed out, the function $\xi\mapsto g_{L,A}(\xi,s)$ is strictly increasing for $s<A$ so we can define its inverse, $h_{L,A}(\xi,s)$, by
\begin{equation}
%\label{}
g_{L,A}\bigl(h_{L,A}(\xi,s),s\bigr)=\xi.
\end{equation}
Then we have:

\begin{lemma}
%\label{lemma}
Let~$f$ be the probability density of~$\xi(0)$. For any event~$G$ depending only on $\{\xi(x)\}_{x\in C}$, any $A>0$, any $s\in[0,A)$ and all $L$ sufficiently large,
\begin{equation}
\label{E:6.10q}
\BbbP\bigl(\tilde\xi_{s}\in G\bigr)
=\E\Biggl(\,\,\1_{G}\,\,\Bigl(\frac A{A-s}\Bigr)^{K_{L,A}(\xi,s)}\prod_{x\in C}\frac{f(h_{L,A}(\xi(x),s))}{f(\xi(x))}\Biggr),
\end{equation}
where $K_{L,A}(\xi,s):=\#\{x\in C\colon A<\hata_L-h_{L,A}(\xi(x),s)<2A\}$.
\end{lemma}

\begin{proofsect}{Proof}
Let $L$ be so large that the probability density~$f$ is well defined and positive for all arguments larger than $\hata_L-2A$. Notice the change of variables $\xi\mapsto\tilde\xi_s$, with explicit form $\xi=h_{L,A}(\tilde\xi_s,s)$, incurs the Radon-Nikodym derivative
\begin{equation}
%\label{}
\frac{\textd\xi}{\textd\tilde \xi_s}=\Bigl(\frac A{A-s}\Bigr)^{\1_{\hata_L-\xi\in(A,2A)}}
\end{equation}
for the corresponding Lebesgue measures. Multiplying this by the ratio of the probability densities, ${f(\xi)}/{f(\tilde\xi_s)}$ gives us the Radon-Nikodym derivative of the law of $\xi$ with respect to the law of~$\tilde\xi_s$.  The result thus follows by writing $\BbbP(\tilde\xi_{s}\in G)$ as an integral with respect to the Lebesgue measure $\prod_{x\in C}\textd\xi(x)$ and changing variables using $\xi(x)=h_{L,A}(\tilde\xi_s(x),s)$ for each~$x$.
\end{proofsect}

We will now proceed to deal with the Radon-Nikodym terms in \eqref{E:6.10q}. The ratios of the probability densities will be controled using  \eqref{densityass} as follows:

\begin{lemma} 
\label{lemma-4.13}
Let~$f$ denote the probability density of~$\xi(0)$. For any finite $C\subset\Z^d$, any $\varphi=(\varphi(x))_{x\in C}$ and $\alpha=(\alpha(x))_{x\in C}$ with $\alpha\ge0$, there exists a quantity $o(1)$ such that
\begin{equation}\label{proddensity}
\prod_{x\in C}\frac{f(\hata_L+\varphi(x)-b_L\alpha(x))}{f(\hata_L+\varphi(x))}
=\exp\biggl\{\bigl(1+o(1)\bigr)\sum_{x\in C}\alpha(x)\texte^{\varphi(x)/\rho}\biggr\}.
\end{equation}
Moreover, $o(1)\to0$ as~$L\to\infty$ uniformly in $C$ and $\varphi,\alpha\in I^{\,C}$, $\alpha\ge0$, for any compact~$I\subset\R$.
\end{lemma}

\begin{proofsect}{Proof} 
For~$F$ be as in Assumption~\ref{ass} and $t\in\R$ such that $\hata_L+t$ lies in the domain of~$F$, let $h_L(t)$ be defined by
\begin{equation}
\label{E:6.10}
\rho\bigl(1+h_L(t)\bigr)\texte^{t/\rho}:=b_L\texte^{F(\hata_L+t)}.
\end{equation}
Thanks to \eqref{E:4.46q}, Assumption~\ref{ass} and the Mean-Value Theorem for~$F$, we have $h_L(t)\to0$ as $L\to\infty$ locally uniformly in~$t$. Next, for~$u\ge0$ such that $f(\hata_L+t-b_Lu)$ and $f(\hata_L+t)$ are well-defined and positive, let $q_L(t,u)$ be defined by
\begin{equation}
\label{E:6.11}
\frac{f(\hata_L+t-b_Lu)}{f(\hata_L+t)} =: \exp\Bigl\{\bigl(1+q_L(t,u)\bigr) u\,\texte^{t/\rho}\Bigr\}.
\end{equation}
Using \eqref{densityass} with $r:=\hata_L+t$ and $s:=b_Lu\,\texte^{F(\hata_L+t)}$, and applying \eqref{E:6.10}, we get $q_L(t,u)\to0$ locally uniformly in $t$ and $u\ge0$. 

Thanks to $u\ge0$, \eqref{E:6.11} can be written as upper/lower bounds valid for~$L$ large once $t,u$ are confined to compact sets in~$\R$. Setting $t:=\varphi(x)$ and $u:=\alpha(x)$ and applying this bound to the product in \eqref{proddensity}, the claim follows.
\end{proofsect}

We are now ready to prove the main claim of this subsection:

\begin{proofsect}{Proof of Proposition~\ref{prop-smallbox}}
Fix $\epsilon>0$ and $A>0$ and set
\begin{multline}
%\label{}
\qquad
G:=\bigl\{\xi\colon\lambda_{B_{R_L}}^{\ssup 1}(\xi)\ge a_L+rb_L\bigr\}
\cap\bigl\{\xi\colon \LL_{B_{R_L},A}(\xi-\hata_L)\ge1-\epsilon\bigr\}
\\
\cap\bigl\{\xi\colon \LL_{B_{R_L},2A}(\xi-\hata_L)\le1+\epsilon\}
\cap\bigl\{\xi\colon\max_{x\in B_{R_L}}\xi(x)\le\hata_L+A\bigr\}.
\qquad
\end{multline}
Our ultimate goal is to show that the right hand side of \eqref{E:6.10q} with this~$G$ and~$s$ replaced by $sb_L$ is asymptotically equal to $\texte^{-s+o(1)}\BbbP(\xi\in G)$.

The first term in the expectation in \eqref{E:6.10q} is bounded directly: Since $s\ge0$ we have, for~$L$ sufficiently large and all $\xi$, that
\begin{equation}
%\label{}
1\le\Bigl(\frac A{A-s}\Bigr)^{K_{L,A}(\xi,s)}\le \texte^{O(1)s b_L R_L^d},
\end{equation}
where we estimated $K_{L,A}(\xi,s)$ by the total volume of~$B_{R_L}$. Since $R_L=o(\log L)$ while $b_L=O(1/\log L)$, the right-hand side tends to one uniformly on compact sets of~$s$.

For the product of ratios of probability densities, we will apply Lemma~\ref{lemma-4.13}. Given a configuration $\xi$, let us abbreviate $\varphi(x):=\xi(x)-A$ and define $\alpha(x)$ by $\xi(x)-sb_L\alpha(x):=h_{L,A}(\xi(x),sb_L)$. As is easy to check, $\alpha(x)\in[0,1]$ while $\varphi(x)\in[-2A,A]$ for all~$\xi\in G$ where $h_L(\xi(x),sb_L)\ne\xi(x)$. Lemma~\ref{lemma-4.13} thus implies
\begin{equation}
%\label{}
\prod_{x\in B_{R_L}}\frac{f(h_{L,A}(\xi(x),sb_L))}{f(\xi(x))}=\exp\Bigl\{\,s\bigl(1+o(1)\bigr)\sum_{x\in B_{R_L}}\!\!\alpha(x)\,\texte^{\varphi(x)/\rho}\Bigr\},
\end{equation}
where $o(1)\to0$ as $L\to\infty$ uniformly on~$G$. Concerning the sum in the exponential on the right, here we note that $\alpha(x)=1$ when $\xi(x)\ge\hata_L-A$ while $\alpha(x)=0$ when $\xi(x)<\hata_L-2A$. Hence,
\begin{equation}
%\label{}
\LL_{B_{R_L},A}(\xi-\hata_L)\le\sum_{x\in B_{R_L}}\alpha(x)\texte^{\varphi(x)/\rho}\le \LL_{B_{R_L},2A}(\xi-\hata_L).
\end{equation}
On~$G$ the left-hand side is at least $1-\epsilon$ while the right-hand side is at most~$1+\epsilon$.
We conclude
\begin{equation}
\label{E:4}
\BbbP\bigl(\tilde\xi_{sb_L}\in G\bigr)=\texte^{-s+O(\epsilon)}\BbbP(\xi\in G),
\end{equation}
where $O(\epsilon)$ is bounded by a constant times~$\epsilon$ uniformly on compact sets of~$s$, for all~$A>0$ sufficiently large (and larger than~$s$).

We are ready to put all the above together and extract the desired claim. First, Lemmas~\ref{lemma-bigcup},~\ref{lemma-4.9} and~\ref{lemma-3.6q} and the bound $\BbbP(\max_{x\in B_{R_L}}\xi(x)>\hata_L+A)=o(L^{-d})$ yield
\begin{equation}
%\label{}
\BbbP(\xi\in G)=\BbbP\bigl(\lambda_{B_{R_L}}^{\ssup 1}(\xi)\ge a_L+rb_L\bigr)+o(L^{-d}),
\qquad L\to\infty.
\end{equation}
Since $\xi\ge\tilde\xi_s\ge\xi-sb_L$, we similarly get
\begin{equation}
%\label{}
\BbbP\bigl(\tilde\xi_{sb_L}\in G\bigr)=
\BbbP\bigl(\lambda_{B_{R_L}}^{\ssup 1}(\tilde\xi_{sb_L})\ge a_L+rb_L\bigr)+o(L^{-d}),
\qquad L\to\infty.
\end{equation}
Plugging these into \eqref{E:4}, invoking the inclusions \twoeqref{E:4.49w}{E:4.50w}  and noting that~$s'$ in \eqref{E:4.50w} can be made arbitrarily close to~$s$ by increasing~$A$, we conclude the claim for~$s\ge0$. For $s<0$ the claim follows by symmetry.
\end{proofsect}

\subsection{Stability with respect to partition size}
In order to conclude the proof of Theorem~\ref{thm-maxorder}, we need to relate the upper tails of the law of the principal eigenvalues in~$B_{R_L}$ and~$B_{N_L}$. Related to this is the question on how much does~$a_L$, defined in \eqref{E:4.13q}, depend on the (rather arbitrary) choice of the sequence~$N_L$. As attested by the next lemma, one direction is quite easy:

\begin{lemma}
\label{lemma-4.13}
There exists a constant $c=c(d)\in(0,\infty)$ such that or any $N\ge R$ and any $a\in\R$,
\begin{equation}
\label{E:7.22eq}
-\log\bigl(1-\BbbP(\lambda^{\ssup 1}_{B_N}\ge a)\bigr)
\ge\bigl(1-cR/N\bigr)\Bigl(\frac NR\Bigr)^d\BbbP\bigl(\lambda^{\ssup 1}_{B_{R}}\ge a\bigr).
\end{equation}\end{lemma}

\begin{proofsect}{Proof}
Let us cover $\Z^d$ by disjoint translates of~$B_R$ and let $B_R^{(i)}$, $i=1,\dots,n$, denote those translates that are contained in the box~$B_N$. Then $\lambda^{\ssup 1}_{B_N}(\xi)\ge\lambda^{\ssup 1}_{B_R^{(i)}}(\xi)$ for every~$i$ and since $\lambda^{\ssup 1}_{B_R^{(i)}}(\xi)$ are independent and equidistributed to $\lambda^{\ssup 1}_{B_R}(\xi)$, we thus have
\begin{equation}
\begin{aligned}
\BbbP\bigl(\lambda^{\ssup 1}_{B_N}\ge a\bigr)
&\ge 1-\BbbP\bigl(\lambda^{\ssup 1}_{B_R}< a\bigr)^n
\\
&\ge 1-\exp\bigl\{-n\,\BbbP(\lambda^{\ssup 1}_{B_R}< a)\bigr\}.
\end{aligned}
\end{equation}
The claim follows by taking a log and using that $n\ge(1-cR/N)(N/R)^d$ for some~$c>0$.
\end{proofsect}

Notice that \eqref{E:7.22eq} implies that once $\BbbP(\lambda^{\ssup 1}_{B_{N_L}}\ge t_L)\to0$ and $R_L/N_L\to0$ as $L\to\infty$, for some sequences~$R_L$, $N_L$ and $t_L$, then also
\begin{equation}
%\label{}
\BbbP(\lambda^{\ssup 1}_{B_{N_L}}\ge t_L)\ge\bigl(1+o(1)\bigr)\Bigl(\frac{N_L}{R_L}\Bigr)^d\BbbP\bigl(\lambda^{\ssup 1}_{B_{R_L}}\ge t_L\bigr),\qquad L\to\infty.
\end{equation}
(Indeed, just expand the log into a power series and dominate it by the first-order term.) 

The bound in the opposite direction will require introducing an auxiliary scale~$R'_L$ as follows: Suppose, for the sake of present section, that $R_L$ and $N_L$ are sequences of integers such that 
\begin{equation}
%\label{}
\lim_{L\to\infty}\,\frac{R_L}{\log\log L}=\infty,\quad \lim_{L\to\infty}\,\frac{R_L}{N_L}=0
\quad\text{and}\quad\lim_{L\to\infty}\,\frac{N_L}{L}=0
\end{equation}
and let $R'_L$ be a sequence of integers satisfying
\begin{equation}
\label{E:6.26}
\lim_{L\to\infty}\,\frac{R'_L}{R_L}=0.
\end{equation}
Then we have:

\begin{lemma}
\label{lemma-6.7}
For any $A>0$ and any sequence $t_L\ge-A$ there is $c>0$ such that
\begin{equation}
%\label{}
\BbbP\bigl(\lambda^{\ssup 1}_{B_{N_L}}\ge \hata_L+t_L\bigr)
\le o(L^{-d})+\bigl(1+o(1)\bigr)\Bigl(\frac{N_L}{R_L}\Bigr)^d\BbbP\bigl(\lambda^{\ssup 1}_{B_{R_L}}\ge \hata_L+t_L-\texte^{-cR'_L}\bigr),
\end{equation}
as $L\to\infty$.
\end{lemma}

\begin{proofsect}{Proof}
Pick~$A>0$ such that $t_L\ge -A$ and consider the set $\mathfrak C$ of connected components of the union of balls~$B_{R'_L}(x)$ for~$x\in B_{N_L}$ such that $\xi(x)\ge\hata_L-3A$. By Lemma~\ref{lemma-4.7a} (with~$R_L$ replaced by~$R'_L$) and $N_L\le L$, there is an integer~$n_A>0$ such that
\begin{equation}
\label{E:6.28}
\BbbP\bigl(\max_{\CC\in\mathfrak C}\,\text{diam}\,\CC>n_A R'_L\bigr)=o(L^{-d}).
\end{equation}
Now consider a partition of~$\Z^d$ into disjoint translates of~$B_{R_L}$ and let~$B_{R_L}^{(i)}$, $i=1,\dots,m_L$, denote those boxes in the covering that have at least one vertex in common with~$B_{N_L}$. Considering the set $S$ of all vertices on the inner boundary of these boxes, let~$B_{2n_AR'_L}^{(j)}$, $j=1,\dots,k_L$, denote a covering thereof by translates of~$B_{2n_AR'_L}$ centered at these vertices such that no vertex in~$S$ lies in more than two boxes from these. The key point is that, on the event
\begin{equation}
%\label{}
G:=\bigl\{\xi\colon\max_{\CC\in\mathfrak C}\,\text{diam}\,\CC\le n_AR'_L\bigr\},
\end{equation}
each component $\CC\in\mathfrak C$ is entirely contained in one of the above boxes $B_{R_L}^{(i)}$ or $B_{2n_AR'_L}^{(j)}$.

Since $\xi(x)\ge\lambda_{B_{N_L}}^{\ssup1}(\xi)-2A$ and $\lambda_{B_{N_L}}^{\ssup1}(\xi)\ge\hata_L+t_L$ imply $\xi(x)\ge\hata_L-3A$,  Theorem~\ref{thm-truncation} can be used for the set $U:=\bigcup_{\CC\in\mathfrak C}\CC$. Thereby we get
\begin{equation}
%\label{}
\lambda_{B_{N_L}}^{\ssup1}(\xi)\le\max_{\CC\in\mathfrak C}\lambda_{\CC}^{\ssup1}(\xi)+\epsilon_{R'_L},
\end{equation}
where $\epsilon_{R'_L}:=2d(1+\frac A{2d})^{1-2R'_L}$. The monotonicity of $C\mapsto\lambda_{C}^{\ssup1}(\xi)$ then shows that
\begin{equation}
\label{E:6.32}
\lambda_{B_{N_L}}^{\ssup1}(\xi)\le\epsilon_{R'_L}+
\max\biggl\{\,\max_{i=1,\dots,m_L}\lambda_{B_{R_L}^{(i)}}^{\ssup1}(\xi),\max_{j=1,\dots,k_L}\lambda_{B_{2n_AR'_L}^{(j)}}^{\ssup1}(\xi)\biggr\}
\end{equation}
holds on $G\cap\{\lambda_{B_{N_L}}^{\ssup1}(\xi)\ge\hata_L+t_L\}$.

Applying the union bound \eqref{E:6.32} and \eqref{E:6.28} yield
\begin{equation}
\begin{aligned}
\BbbP\bigl(\lambda^{\ssup 1}_{B_{N_L}}\ge \hata_L+t_L\bigr)
&\le m_L\,\BbbP\bigl(\lambda^{\ssup 1}_{B_{R_L}}\ge \hata_L+t_L-\epsilon_{R'_L}\bigr)
\\
&\qquad\qquad+k_L\,\BbbP\bigl(\lambda^{\ssup 1}_{B_{2n_AR'_L}}\ge \hata_L+t_L-\epsilon_{R'_L}\bigr)
+o(L^{-d})
\\
&\le (m_L+k_L)\,\BbbP\bigl(\lambda^{\ssup 1}_{B_{R_L}}\ge \hata_L+t_L-\epsilon_{R'_L}\bigr)+o(L^{-d}),
\end{aligned}
\end{equation}
where the last inequality holds because $2n_AR'_L\le R_L$. Since $R'_L\ll R_L$,
\begin{equation}
%\label{}
m_L=\bigl(1+o(1)\bigr)\Bigl(\frac{N_L}{R_L}\Bigr)^d\quad\text{and}\quad
k_L\le O(1)\frac{N_L}{R_L}\,\Bigl(\frac{N_L}{R'_L}\Bigr)^{d-1}=o(m_L)
\end{equation}
the claim follows by noting that $\epsilon_{R'_L}\le\texte^{-cR'_L}$ for some $c>0$.
\end{proofsect}

\subsection{Proof of eigenvalue order statistics}
First we establish the Gumbel max-order tail for the principal eigenvalues:

\begin{proofsect}{Proof of Theorem~\ref{thm-maxorder}}
Choose~$r_L$ so that $r_L/\log\log L\to\infty$ and \eqref{E:6.26} hold. Then $\texte^{-cr_L}=o(b_L)$ as~$L\to\infty$, for any~$c>0$. Combining Lemma~\ref{lemma-6.7} and Proposition~\ref{prop-smallbox}, the claim follows.
\end{proofsect}

This then implies the extreme order law for eigenvalues:

\begin{proofsect}{Proof of Corollary~\ref{cor-1.3}}
This is a direct consequence of Theorems~\ref{prop-4.6} and~\ref{thm-maxorder}, the standard results about max-order statistics and the facts that $R_L/\log\log L\to\infty$ and $\log|D_L|=(d+o(1))\log L$ as~$L\to\infty$ for any $D\in\mathfrak D$.
\end{proofsect}

\section{Eigenfunction decay}
\label{sec7}\noindent
In this short section we will provide the arguments need in the proof of eigenfunction localization. Recall our notation $\epsilon_R:=2d(1+\frac A{2d})^{1-2R}$. A key observation is:

\begin{lemma}
\label{lemma-7.1}
Let~$R_L/\log\log L\to\infty$. Then for each $A>0$ and each $k\ge2$,
\begin{equation}
%\label{}
\frac1{\epsilon_{R_L}}\,
\min_{\ell=1,\dots,k-1}\bigl(\lambda_{D_L}^{\ssup\ell}(\xi)-\lambda_{D_L}^{\ssup{\ell+1}}(\xi)\bigr)
\,\underset{L\to\infty}\longrightarrow\infty
\end{equation}
in probability.
\end{lemma}

\begin{proofsect}{Proof}
By the convergence to the Poisson point process established in Corollary~\ref{cor-1.3},
\begin{equation}
%\label{}
\bigl[\lambda_{D_L}^{\ssup\ell}(\xi)-\lambda_{D_L}^{\ssup{\ell+1}}(\xi)\bigr]\log L\,\,\underset{L\to\infty}\Longrightarrow\,\,\frac\rho d\,\log\frac{Z_1+\dots+Z_{\ell+1}}{Z_1+\dots+Z_\ell},
\end{equation}
where $Z_1,Z_2,\dots$ are i.i.d.\ exponentials with parameter one. Since the right-hand side is positive with probability one, and $\epsilon_{R_L}\log L\to0$, the result follows.
\end{proofsect}

Next let us consider the distance $\textd(x,\CC)$ defined, as an example, right before Theorem~\ref{prop-exp-decay} and let $\dist(x,y)$ stand for the $\ell^1$-distance between~$x$ and~$y$ on~$\Z^d$. Clearly, $\textd(x,\CC)\le\dist(x,\CC)$. Our next item of business is to show:

\begin{lemma}[Comparison of distances]
\label{lemma-7.2}
For each~$A>0$ there are $c_0,c_1,c_2\in(0,\infty)$ such that for any $R_L\to\infty$ that satisfies $R_L\le c_0\log L$, we have
\begin{equation}
%\label{}
\textd(x,\CC)\ge c_1\,\text{\rm dist}(x,\CC)-c_2 R_L,\qquad x\in B_L,\,\CC\in\mathfrak C_{R_L,A},
\end{equation}
holds with probability tending to one as $L\to\infty$.
\end{lemma}

\begin{proofsect}{Proof}
By its definition, $\textd(x,y)$ is the (graph-theoretical) distance on the graph obtained by contracting each connected components~$\CC$ to a single vertex. So it suffices to prove 
\begin{equation}
%\label{}
\textd(x,y)\ge c_1\,\text{\rm dist}(x,y)-c_2 R_L,\qquad x,y\in B_L.
\end{equation}
Let~$G:=\{\xi\colon\diam\CC\le n_A R_L\,\,\forall\CC\}$, where the diameter is in the $\ell^1$-distance on~$\Z^d$.
Given~$x$ and~$y$, consider a path~$\pi$ on this contracted graph achieving~$\textd(x,y)$. This can be extended into a path on~$B_L$ by concatenating with paths inside the components, which yields
\begin{equation}
%\label{}
\dist(x,y)\le\textd(x,y)+Yn_A R_L,\qquad \text{on }G,
\end{equation}
where $Y$ denotes the number of connected connected components encountered by~$\pi$.

To estimate~$Y$, consider any vertex self-avoiding path from~$x$ to~$y$ and let~$\KK$ denote the union of $B_{R_L}(z)$ for all~$z$ on this path. Clearly, $|\KK|\le c \dist(x,y) R_L^d$ for some constants~$c>0$. By a union bound and \eqref{5.1-eq},
\begin{equation}
%\label{}
\BbbP\bigl(Y\ge n)\le |\KK|L^{-n\theta}\le c \dist(x,y) R_L^d L^{-nd\theta}\le L^{-n\theta'}
\end{equation}
for some $\theta,\theta'>0$ and all $n\ge1$, where we used that $\dist(x,y) R_L^d\le L^{1+o(1)}$. Hence, for any $\eta>0$,
\begin{equation}
%\label{}
\BbbP\Bigl(Y\ge \eta\frac{\dist(x,y)}{R_L}\Bigr)\le \exp\Bigl\{-\theta'\eta\frac{\log L}{R_L}\dist(x,y)\Bigr\}.
\end{equation}
Summing this over $x,y\in B_L$ with $\dist(x,y)\ge R_L/\eta^2$, the result will tend to zero with~$L\to\infty$ provided~$\eta$ is sufficiently small. As also $G$ has probability tending to one,  we get
\begin{equation}
%\label{}
\dist(x,y)\le\textd(x,y)+\eta\frac{\dist(x,y)}{R_L}n_A R_L,
\end{equation}
implying $\textd(x,y)\ge (1-\eta n_A)\dist(x,y)$ as soon as $\dist(x,y)>R_L/\eta^2$, with probability tending to one as $L\to\infty$. As $\eta$ can be chosen so that $\eta n_A<1$, we are done.
\end{proofsect}

We are now ready to establish the eigenfunction decay starting first with long distances:

\begin{proofsect}{Proof of Theorem~\ref{thm2.4}, large distances}
We will prove \eqref{E:1.20q} at distances at least~$\log L/\log\log L$. Our aim is to apply Theorem~\ref{prop-exp-decay} for $\lambda:=\lambda_{D_L}^{\ssup k}(\xi)$, any $A>0$ and~$R:=R_L$, where
\begin{equation}
%\label{}
R_L:=\lfloor\log L/\log\log L\rfloor.
\end{equation}
We will now check the conditions (1-3) of Theorem~\ref{prop-exp-decay}.

First, since $\epsilon_{R_L}=o(1/\log L)$, condition~(1) holds thanks to Lemma~\ref{lemma-7.1}. Concerning condition~(2), we note that as soon as
\begin{equation}
\label{E:7.9a}
R_L\,\frac{\log\bigl(\BbbP(\xi(0)>t_L)^{-1}\bigr)}{\log L}\,\underset{L\to\infty}\longrightarrow\,\,\infty
\end{equation}
for some sequence $t_L$, then the following holds with probability tending to one as $L\to\infty$: For any self-avoiding path $(x_1,\dots,x_n)$ in~$D_L$ of length $n\ge R_L$,
\begin{equation}
%\label{}
\#\bigl\{i=1,\dots,n\colon \xi(x_i)\ge t_L\bigr\}\le\frac n2.
\end{equation}
Assumption~\ref{ass} tells us that $\log(\BbbP(\xi(0)>a)^{-1})=\texte^{F(a)}$ with $F(a)=(1/\rho+o(1))a$ as~$a\to\infty$, we easily check that \eqref{E:7.9a} holds for, say, $t_L:=\hata_L/2$. Since $\lambda=\hata_L-o(1)$, the condition in \eqref{E:4.40} is valid with~$h:=c\log\log L$ for some~$c>0$. 

Concerning~(3), by Lemma~\ref{lemma-4.7a} (which holds as soon as~$R_L=o(L)$) all components of $\mathfrak C_{R_L,A}$ have diameter at most~$n_A R_L$, with probability tending to one as~$L\to\infty$. Condition~(3) is then readily checked for any~$\delta\in(0,1)$.

Since the premises of Theorem~\ref{prop-exp-decay} hold, we know that there is a component $\CC\in\mathfrak C_{R_L,A}$ such that $\psi:=\psi_{D_L,\xi}^{\ssup k}$ obeys
\begin{equation}
%\label{}
\bigl|\psi(z)\bigr|\le\texte^{-c\delta(\log\log L)\textd(z,\CC)},\qquad z\in D_L.
\end{equation}
In particular, $X_k$, defined by \eqref{E:1.9w}, must satisfy~$\dist(X_k,\CC)\le R_L$. As Lemma~\ref{lemma-7.2} is at our disposal, we further conclude
\begin{equation}
\begin{aligned}
\textd(z,\CC)
&\ge c_1\dist(z,\CC)-c_2 R_L 
\\
&\ge c_1|z-X_k|-(c_1n_A+c_2+1)R_L,
\end{aligned}
\end{equation}
where we used that, by Lemma~\ref{lemma-4.7a}, $\diam\CC\le n_A R_L$ with probability tending to one as~$L\to\infty$. Hereby we get
\begin{equation}
\label{E:7.13a}
|z-X_k|\ge\frac{c_1n_A+c_2+1}{2c_1}R_L\quad\Rightarrow\quad
\bigl|\psi(z)\bigr|\le\texte^{-c'(\log\log L)\dist(z,\CC)},
\end{equation}
for~$c'$ given by~$c':=c\delta c_1/2$. In particular, \eqref{E:1.20q} is true.
\end{proofsect}

Before we move on to the short distances, let us abbreviate $R_L':=\frac{c_1n_A+c_2+1}{2c_1}R_L$ for~$R_L$ as in the previous proof and notice that \eqref{E:7.13a} yields
\begin{equation}
\sum_{|z-X_k|\ge R_L'}\xi(z)\bigl|\psi(z)\bigr|^2=o(1),
\end{equation}
with probability tending to one as~$L\to\infty$ because $\xi(z)$ is at most a constant times $\log\log L$ in this limit. In particular, we have
\begin{equation}
\label{E:7.14a}
\lambda_{B_{R_L'}(X_k)}^{\ssup1}(\xi)=\lambda-o(1)
\end{equation}
by the variational characterization of the principal eigenvalue.

\begin{proofsect}{Proof of Theorem~\ref{thm2.4}, short distances}
It remains to prove exponential decay for distances less than order~$R_L$. Set $R_L'$ as right before this proof. Our aim is to use Lemma~\ref{lemma-3.8w} to show that all but a finite number of values in~$B_{R_L'}(X_k)$ are more than a positive constant below~$\lambda$. We will now proceed to verify the premises of Lemma~\ref{lemma-3.8w}. 

Pick~$\epsilon>0$. Thanks to Lemmas~\ref{lemma-bigcup} and~\ref{lemma-4.9}, and $\chi_{B_{R_L}}\downarrow\chi$, the following holds with probability one as~$L\to\infty$: Once~$A$ is sufficiently large,
\begin{equation}
%\label{}
\LL_{C,A}\bigl(\xi-(a_L-\epsilon)-\chi_C\bigr)\le\texte^{\epsilon/(2\rho)}
\end{equation}
holds for~$C$ ranging over all translates of~$B_{R_L'}$ that intersect~$D_L$. Set $\delta:=\epsilon/(2\rho)$ and abbreviate $C:=B_{R_L'}(X_k)$. Next observe that for $A':=-\frac12\rho\log(2\sinh\delta)$ and~$\epsilon$ small, there is~$r\ge1$ so that
\begin{equation}
%\label{}
2d\Bigl(1+\frac{A'-d}{2d}\Bigr)^{1-2r}\le\frac\epsilon2.
\end{equation}
(Indeed, $A'=\frac12\log(2\delta)+O(1)$ so setting $r$ to be proportional to a constant times $\log(1/\epsilon)$ will do once~$\epsilon$ is small enough.) As $\lambda=a_L+o(1)$, the bound \eqref{E:7.14a} tells us 
\begin{equation}
%\label{}
\lambda_{C}^{\ssup1}(\xi)\ge a_L-\frac\epsilon2\ge a_L-\epsilon+2d\Bigl(1+\frac{A'-d}{2d}\Bigr)^{1-2r}.
\end{equation}
By taking~$\delta$ small and~$A$ large, and setting $A'':=A-d$ (which is less than~$A-\chi_C$) the premises of Lemma~\ref{lemma-3.8w} are thus satisfied and so we conclude that $\xi(z)\le\lambda_C^{\ssup1}(\varphi)-A''$ for all~$z\in C$ that are at least~$cr$ away from~$X_k$. As~$A''>0$, Lemma~\ref{lemma-2-norm} then shows
\begin{equation}
%\label{}
\sum_{|z-X_k|\ge a}\bigl|\psi(z)\bigr|^2\le \Bigl(1+\frac{A''}{2d}\Bigr)^{-2a}\Vert\psi\Vert_2^2
\end{equation}
for all~$a>cr$. In particular, $\psi(z)$ decays exponentially with distance from~$X_k$.
\end{proofsect}

Finally, we will supply a formal proof of our principal result:

\begin{proofsect}{Proof of Theorem~\ref{thm-orderstat}}
Part~(1) is a direct consequence of the bounds in Theorem~\ref{thm2.4}. For part~(2) notice that argument producing the coupling to i.i.d.\ random variables in Theorem~\ref{prop-4.6} is such that $\lambda_{D_L}^{\ssup k}$ is coupled exactly to $\lambda_{B_{N_L}^{(i)}}(\xi)$ for~$i$ such that $X_k\in B_{N_L}^{(i)}$. In the notation of this theorem, since $N_L=o(L)$, the collection \eqref{E:1.8} is well approximated by
\begin{equation}
%\label{}
\biggl\{\Bigl(\frac{z_i}L,\frac1\rho\bigl(\lambda_i(\xi)-a_L\bigr)\log|D_L|\Bigr)\colon i=1,\dots,m_L\biggr\},
\end{equation}
where $z_i$ denotes the point at the center of $B_{N_L}^{(i)}$. As~$m_L\to\infty$ while the number of points with the second coordinate above a given value is with high probability bounded, this is in turn well approximated by sampling the first coordinate uniformly from $\{z_i\colon i=1,\dots,m_L\}$, independently of the second coordinate. The last approximating process converges to a Poisson point process with intensity measure $\1_D\textd x\otimes\texte^{-\lambda}\textd\lambda$.
\end{proofsect}

\section*{Acknowledgments}
\noindent
This research has been partially supported by DFG Forschergruppe 718 ``Analysis and Stochastics in Complex Physical Systems,'' NSF grants DMS-0949250 and DMS-1106850, NSA grant H98230-11-1-0171 and GA\v CR project P201-11-1558. The final parts of this paper were written when M.B.\ was visiting RIMS at Kyoto University, whose hospitality is gratefully acknowledged.

\renewcommand{\sc}{}
%%%%%%%%%%%%%%%%%%%%%%%%%%%%%%%%%%%%%%%%%%%%%%%%

\end{document}